\documentclass{gtart}

\def\ifplaintex{\expandafter\ifx\csname documentclass\endcsname\relax}

\def\gtp{{\mathsurround=0pt\it $\cal G\mskip-2mu$eometry \&\ 
$\cal T\!\!$opology $\cal P\!$ublications}}  

\def\recd{{\small Received:\qua\receiveddate\ifx\reviseddate\relax
\else\qquad Revised:\qua\reviseddate\fi\par}} 

\def\lognumber#1{\def\thelognumber{#1}}
\def\volumenumber#1{\def\thevolumenumber{#1}}
\def\volumeyear#1{\def\thevolumeyear{#1}}
\def\papernumber#1{\def\thepapernumber{#1}}
\def\pagenumbers#1#2{\def\startpage{#1}\def\finishpage{#2}}
\def\published#1{\def\publishdate{#1}}

\def\received#1{\def\receiveddate{#1}}
\def\revised#1{\def\reviseddate{#1}}
\def\accepted#1{\def\accepteddate{#1}}

\long\def\asciiabstract#1{\long\def\theasciiabstract{#1}}

\let\\\par\let\thelognumber\relax\let\thevolumenumber\relax
\let\thepapernumber\relax\let\thevolumeyear\relax\let\startpage\relax
\let\finishpage\relax\let\publishdate\relax\let\receiveddate\relax
\let\reviseddate\relax\let\accepteddate\relax\let\theasciititle\relax
\let\theasciiauthors\relax
\let\theasciiabstract\relax

\let\theasciiemail\relax

\ifplaintex
\font\logobig=cmssbx10 scaled 3836
\font\logomed=cmssbx10 scaled 2557
\else
\font\logobig=cmssbx10 scaled 4200
\font\logomed=cmssbx10 scaled 2800
\fi

\long\def\makeagttitle{   
\count0=\startpage
\agt\hfill      
\hbox to 45truept{\vbox to 0pt{\vglue -13truept{\logomed A\kern -.37em{\logobig 
T}\kern -.38em G}\vss}\hss}
\break
{\small Volume \thevolumenumber\ (\thevolumeyear)
\startpage--\finishpage\nl
Published: \publishdate}

\vglue .25truein

{\parskip=0pt\leftskip 0pt plus
1fil\def\\{\par\smallskip}{\Large\bf\thetitle}\par\medskip} \vglue
0.05truein

{\parskip=0pt\leftskip 0pt plus 1fil\def\\{\par}{\sc\theauthors}
\par\medskip}%
 
\vglue 0.03truein 

{\small\leftskip 25truept\rightskip 25truept{\bf Abstract}\stdspace\theabstract

{\bf AMS Classification}\stdspace\theprimaryclass
\ifx\thesecondaryclass\relax\else; \thesecondaryclass\fi\par
{\bf Keywords}\stdspace \thekeywords\par}\vglue 7truept

}   

\ifplaintex
\font\phead=cmsl9 scaled 950
\font\pnum=cmbx10 scaled 913
\font\pfoot=cmsl9 scaled 950
\headline{\vbox to 0pt{\vskip -4.5mm\line{\small\phead\ifnum
\count0=\startpage ISSN 1472-2739 (on-line) 1472-2747 (printed)
\hfill {\pnum\folio}\else\ifodd\count0\def\\{ }%
\ifx\theshorttitle\relax\thetitle\else\theshorttitle\fi\hfill{\pnum\folio}
\else\def\\{ and }{\pnum\folio}\hfill\ifx\theshortauthors\relax\theauthors
\else\theshortauthors\fi\fi\fi}\vss}}
\footline{\vbox to 0pt{\vglue 0mm\line{\small\pfoot\ifnum\count0=\startpage
\copyright\ \gtp\hfill\else
\agt, Volume \thevolumenumber\ (\thevolumeyear)\hfill\fi}\vss}}
\else
\font=cmsl9 scaled 1050
\font\lnum=cmbx10 
\font=cmsl9 scaled 1050
\makeatletter
\def\@oddhead{{\small\lhead\ifnum\count0=\startpage ISSN 1472-2739 
(on-line) 1472-2747 (printed)\hfill {\lnum\number\count0}\else\ifodd\count0
\def\\{ }\ifx\theshorttitle\relax \thetitle \else\theshorttitle\fi\hfill
{\lnum\number\count0}\else\def\\{ and }{\lnum\number\count0}
\hfill\ifx\theshortauthors\relax 
\theauthors\else\theshortauthors\fi\fi\fi}}\def\@evenhead{\@oddhead}
\def\@oddfoot{\small\lfoot\ifnum\count0=\startpage\copyright\ \gtp\hfill\else
\agt, Volume \thevolumenumber\ (\thevolumeyear)\hfill\fi}
\def\@evenfoot{\@oddfoot}
\makeatother
\fi
\let\maketitlepage\makeagttitle
\let\makeshorttitle\maketitlepage
\let\maketitle\maketitlepage

\newwrite\gtoutfile
\long\gdef\makeheadfile{  
{\def\\{, }\def\s{ }
\immediate\openout\gtoutfile head.xxx
\immediate\write\gtoutfile{To: math@arxiv.org}
\immediate\write\gtoutfile{Subject: put OR rep NNNNN:ppppp}
\immediate\write\gtoutfile{--text follows this line--}
\immediate\write\gtoutfile{Proxy-for: \ifx\theasciiauthors\relax
\theauthors\else\theasciiauthors\fi\s<\ifx\theasciiemail\relax\theemail\else\theasciiemail\fi>}
\immediate\write\gtoutfile{\noexpand\\}
\immediate\write\gtoutfile{Authors: \ifx\theasciiauthors\relax
\theauthors\else\theasciiauthors\fi}
{\def\\{ }\immediate\write\gtoutfile{Title: \ifx\theasciititle\relax
\thetitle\else\theasciititle\fi}}
\immediate\write\gtoutfile{Subj-class: GT or SG, GR etc}
\immediate\write\gtoutfile{MSC-class: \theprimaryclass\ifx\thesecondaryclass\relax\else, \thesecondaryclass\fi}
\immediate\write\gtoutfile{Journal-ref: Algebr. Geom. Topol. \thevolumenumber\s
(\thevolumeyear) \startpage-\finishpage}
\immediate\write\gtoutfile{Comments: Published by Algebraic and
Geometric Topology at}
\immediate\write\gtoutfile{\s\s\s  http://www.maths.warwick.ac.uk/agt/AGTVol\thevolumenumber/agt-\thevolumenumber-\thepapernumber.abs.html}
\immediate\write\gtoutfile{\noexpand\\}
\immediate\write\gtoutfile{}
\ifx\theasciiabstract\relax
\immediate\write\gtoutfile{\theabstract}\else
\immediate\write\gtoutfile{\theasciiabstract}\fi
\immediate\write\gtoutfile{}
\immediate\write\gtoutfile{\noexpand\\}
\immediate\write\gtoutfile{}
\immediate\closeout\gtoutfile}}  

\def\maketitlepage{\makeagttitle\makeheadfile}
\let\makeshorttitle\maketitlepage
\let\maketitle\maketitlepage

\lognumber{29}
\volumenumber{3}
\volumeyear{2003}
\papernumber{29}
\published{25 September 2003}
\pagenumbers{873}{904}
\received{21 February 2003}
\revised{16 September 2003}
\accepted{23 September 2003}

\usepackage{amsmath,amssymb}
\usepackage{graphics}
\usepackage{psfrag}
\usepackage{amscd}

\newcommand{\Exp}[2]{\ensuremath{\exp_{#1}\!{#2}}}
\newcommand{\s}{\ensuremath{S^1}}
\newcommand{\exps}[1]{\Exp{#1}{\s}}

\newcommand{\expg}[2]{\Exp{#1}{\Gamma_{#2}}}
\newcommand{\expgb}[2]{\Exp{#1}{(\Gamma_{#2},v)}}
\newcommand{\expinf}[1]{\Exp{}{\, {#1}}}
\newcommand{\sym}[2]{\ensuremath{\operatorname{Sym}^{#1}{#2}}}

\newcommand{\pairing}[2]{\langle #1,#2\rangle}

\newcommand{\perm}[1]{\bar{#1}}
\newcommand{\sbinom}[2]{{\genfrac{[}{]}{0pt}{}{#1}{#2}}_{-1}}
\newcommand{\qbinom}[2]{\genfrac{[}{]}{0pt}{}{#1}{#2}}

\newcommand{\drop}{\partial}
\newcommand{\lift}{\delta}
\newcommand{\vcell}[2][{}]{\ensuremath{\sigma^{#2}_{#1}}}
\newcommand{\ecell}[2][{}]{\ensuremath{\tilde\sigma_{#1}^{#2}}}
\def\co{\colon\thinspace}

\newcommand{\J}{\mathcal{J}}        
\newcommand{\C}{\mathcal{C}}        
\newcommand{\cu}{\mathsf{C}}        
\newcommand{\Z}{\mathcal{Z}}        
\newcommand{\V}{\mathsf{V}}         
\newcommand{\F}{\mathcal{F}}        

\newcommand{\real}{\ensuremath{\mathbf{R}}}
\newcommand{\integer}{\ensuremath{\mathbf{Z}}}
\newcommand{\rational}{\ensuremath{\mathbf{Q}}}

\DeclareMathOperator{\intr}{int}
\DeclareMathOperator{\supp}{supp}
\DeclareMathOperator{\Span}{span}
\DeclareMathOperator{\sign}{sign}

\newtheorem{theorem}{Theorem}
\newtheorem{lemma}{Lemma}
\newtheorem{corollary}{Corollary}
\theoremstyle{remark}

\numberwithin{equation}{section}

\begin{document}

\title{Finite subset spaces of graphs\\and punctured surfaces}
\author{Christopher Tuffley}
\address{Department of Mathematics, University of California at Davis \\
         One Shields Avenue, Davis, CA 95616-8633, USA}
\email{tuffley@math.ucdavis.edu}

\begin{abstract}
The $k$th finite subset space of a topological space $X$ is
the space \Exp{k}{X}\ of non-empty finite subsets of $X$ of size at most 
$k$, topologised as
a quotient of $X^k$. The construction is a homotopy functor and
may be regarded as a union of configuration spaces of distinct unordered points
in $X$. We calculate the homology of the finite subset spaces of a 
connected graph $\Gamma$,
and study the maps $(\Exp{k}{\phi})_*$ induced by a map 
$\phi\co\Gamma\rightarrow\Gamma'$ between two such graphs. By homotopy 
functoriality the results apply
to punctured surfaces also. The braid group $B_n$ may be regarded as the
mapping class group of an $n$--punctured disc $D_n$, and as such it acts on
$H_*(\Exp{k}{D_n})$. We prove a structure theorem for this action, showing
that the image of the pure braid group is nilpotent of class at most
$\lfloor (n-1)/2\rfloor$. 
\end{abstract}

\asciiabstract{
The k-th finite subset space of a topological space X is the space
exp_k(X) of non-empty finite subsets of X of size at most k,
topologised as a quotient of X^k. The construction is a homotopy
functor and may be regarded as a union of configuration spaces of
distinct unordered points in X.  We calculate the homology of the
finite subset spaces of a connected graph Gamma, and study the maps
(exp_k(phi))_* induced by a map phi: Gamma --> Gamma' between two such
graphs.  By homotopy functoriality the results apply to punctured
surfaces also. The braid group B_n may be regarded as the mapping
class group of an n-punctured disc D_n, and as such it acts on
H_*(exp_k(D_n)).  We prove a structure theorem for this action,
showing that the image of the pure braid group is nilpotent of class
at most floor((n-1)/2).}


\maketitle

\section{Introduction}

\subsection{Finite subset spaces}

Let $X$ be a topological space and $k$ a positive integer.
The $k$th finite
subset space of $X$ is the space $\Exp{k}{X}$ of nonempty subsets of $X$
of size at most $k$, topologised as a quotient of $X^k$ via the map
\[
(x_1,\ldots,x_k)\mapsto\{x_1\}\cup\cdots\cup\{x_k\}.
\]
The construction is functorial: given a map $f\co X\rightarrow Y$ we obtain
a map $\Exp{k}{f}\co\Exp{k}{X}\rightarrow\Exp{k}{Y}$ by sending 
$S\subseteq X$ to $f(S)\subseteq Y$. Moreover, if $\{h_t\}$ 
is a homotopy between $f$ and $g$ then
$\{\Exp{k}{h_t}\}$ is a homotopy between $\Exp{k}{f}$ and  $\Exp{k}{g}$, 
so that \Exp{k}{}\ is in fact a homotopy functor.

The first finite subset space is of course simply $X$, and
the second finite subset space co-incides with the second symmetric
product $\sym{2}{X}=X^2/S_2$. However, for $k\geq 3$ we have a proper
quotient of the symmetric product as \Exp{k}{X}\ is unable to record
multiplicities: both $(a,a,b)$ and $(a,b,b)$ in $X^3$ 
map to $\{a,b\}$ in \Exp{3}{X}. As a result there are natural
inclusion maps
\begin{equation}
\Exp{j}{X}\hookrightarrow\Exp{k}{X}\co S\mapsto S
\label{include.eq}
\end{equation}
whenever $j\leq k$, stratifying \Exp{k}{X}. 
We define the full finite subset space 
\expinf{X}\ to be the direct limit of this system of inclusions,
\[
\expinf{X} = \varinjlim\Exp{k}{X}.
\]
If $X$ is Hausdorff
then the subspace topology on $\Exp{j}{X}\subseteq\Exp{k}{X}$ co-incides
with the quotient topology it receives from $X^j$~\cite{handel00}.
In this case each stratum $\Exp{j}{X}\setminus\Exp{j-1}{X}$ is 
homeomorphic to the configuration space of sets of $j$ distinct unordered
points in $X$, so that $\Exp{k}{X}$ may be regarded as a union over 
$1\leq j\leq k$ of these spaces. Moreover \Exp{k}{X}\ is compact
whenever $X$ is, in which case it gives a compactification of the 
corresponding configuration space. Such spaces and their compactifications 
have been of considerable interest
recently in algebraic topology: see, for example, Fulton and 
MacPherson~\cite{fulton-macpherson94} and Ulyanov~\cite{ulyanov02}.

For each $k$ and $\ell$ the isomorphism 
$X^k\times X^\ell\rightarrow X^{k+\ell}$ descends to a map
\[
\cup\co\Exp{k}{X}\times\Exp{\ell}{X}\rightarrow\Exp{k+\ell}{X}
\]
sending $(S,T)$ to $S\cup T$. This leads to a form of product on
maps $g\co Y\rightarrow \Exp{k}{X}$, $h\co Z\rightarrow \Exp{\ell}{X}$, and we
define $g\cup h\co Y\times Z\rightarrow \Exp{k+\ell}{X}$ to be the composition
\[
Y\times Z \xrightarrow{g\times h} \Exp{k}{X}\times\Exp{\ell}{X} 
\xrightarrow{\cup} \Exp{k+\ell}{X}.
\]
Clearly $(f\cup g)\cup h = f\cup (g\cup h)$. 
Given a point $x_0\in X$ we obtain as a special case a map
\[
\cup\{x_0\}\co\Exp{k}{X} \rightarrow \Exp{k+1}{X}
\]
taking $S\subseteq X$ to $S\cup\{x_0\}$. The image of $\cup\{x_0\}$
is the subspace $\Exp{k+1}{(X,x_0)}$ consisting of the $k+1$ or fewer
element subsets of $X$ that contain $x_0$. In contrast to the
symmetric product, where the analogous map plays the role 
of~\eqref{include.eq}, the spaces \Exp{k}X\ and \Exp{k+1}{(X,x_0)}\ are in 
general topologically different.
The map $\cup\{x_0\}$ is one-to-one at the point $\{x_0\}\in\Exp{1}{X}$ and on 
the top level stratum $\Exp{k}{X}\setminus\Exp{k-1}{X}$, but is two-to-one 
elsewhere, as $S$ and $S\cup\{x_0\}$ have the same image for $|S|<k$, 
$x_0\not\in S$. Nevertheless the based finite subset spaces
$\Exp{k}{(X,x_0)}$ frequently act as a stepping stone in understanding 
\Exp{k}{X}, often being topologically simpler.

\subsection{History}

The space \Exp{k}{X}\ was introduced by Borsuk and 
Ulam~\cite{borsukulam31} in 1931 as the symmetric product, 
and since then appears to have been
studied at irregular intervals, under various notations, and principally
from the 
perspective of general topology. 
In their original paper Borsuk and Ulam showed that $\Exp{k}{I}\cong I^k$ for
$k=1,2,3$, but that \Exp{k}{I}\ cannot be embedded in $\real^k$ for $k\geq4$.
In 1957 Molski~\cite{molski57} proved similar results for $I^2$ and $I^n$, 
namely that $\Exp{2}{I^2}\cong I^4$ but that neither \Exp{k}{I^2}\ nor 
\Exp{2}{I^k}\ can be embedded in $\real^{2k}$ for any $k\geq 3$. The last was
done by showing that $\Exp{2}{I^k}$ contains a copy of 
$S^k\times\real P^{k-1}$. 

Other authors including Curtis~\cite{curtis86}, Curtis and 
To Nhu~\cite{curtis-nhu85}, Handel~\cite{handel00}, Illanes~\cite{illanes85} 
and Mac\'{\i}as~\cite{macias99} have established general topological and
homotopy-theoretic properties of \Exp{k}{X} and \expinf{X}, and Beilinson and 
Drinfeld~\cite[sec.\ 3.5.1]{chiral} and Ran~\cite{ran93} have used these
spaces in the
context of mathematical physics and algebraic geometry. 
The set \expinf{X}\ has also been studied extensively under a different 
topology as the Pixley-Roy hyperspace of finite subsets of $X$; the
two topologies are surveyed in Bell~\cite{bell79}.
We mention some results on \Exp{k}{X}\ of a homotopy-theoretic
nature. In 1999 Mac\'{\i}as showed that for compact connected metric
$X$ the first \v{C}ech cohomology group $\check{H}^1(\Exp{k}{X};\integer)$
vanishes for $k\geq 3$, and in 2000 Handel proved that for closed
connected $n$--manifolds, $n\geq 2$, the singular cohomology group
$H^i(\Exp{k}{M^n};\integer/2\integer)$ is isomorphic to
$\integer/2\integer$ for $i=nk$, and $0$ for $i>nk$. In addition,
Handel showed that the inclusion maps
$\Exp{k}{(X,x_0)}\hookrightarrow\Exp{2k-1}{(X,x_0)}$ and
$\Exp{k}{X}\hookrightarrow\Exp{2k+1}{X}$ induce the zero map on all
homotopy groups for path connected Hausdorff $X$.

However, although these and
other properties of \Exp{k}{}\ have been established, it appears that 
until recently the only homotopically non-trivial space for which \Exp{k}{X}\ 
was at all well understood for $k\geq 3$ was the circle.
In 1952 Bott~\cite{bott52} proved the surprising result that 
\exps{3}\ is homeomorphic to the three-sphere, correcting Borsuk's
1949 paper~\cite{borsuk49}, and 
Shchepin (unpublished; for three 
different proofs see~\cite{mostovoy99} and~\cite{circles02}) 
later proved the
even more striking result that \exps{1}\ inside \exps{3}\ is a trefoil
knot. An elegant geometric construction due to Mostovoy~\cite{mostovoy99}
in 1999 connects both of these results with known facts about the 
space of lattices in the plane, and in our
previous paper~\cite{circles02} we showed that Bott's and Shchepin's 
results can be viewed as part of a larger pattern: \exps{k}\ has the 
homotopy type of an odd dimensional sphere, and $\exps{k}\setminus\exps{k-2}$
that of a $(k-1,k)$--torus knot complement.
This paper aims to increase the list of spaces for which \Exp{k}{}\ is
understood by using the techniques of~\cite{circles02}\ to study the 
finite subset spaces of connected graphs. The results apply to punctured
surfaces too, by homotopy equivalence, and represent a step towards
understanding finite subset spaces of closed surfaces, as they may be 
used to study these via Mayer-Vietoris type arguments. Further steps
towards this goal are taken in our dissertation~\cite{dissertation}, in
which this paper also appears.

Various different notations have been used for $\Exp{k}{X}$, including 
$X(k)$, $X^{(k)}$, $\mathcal{F}_k(X)$ and $Sub(X,k)$.
Our notation follows that used by 
Mostovoy~\cite{mostovoy99} and reflects the idea that we are truncating
the (suitably interpreted) series
\[
\exp X = \emptyset \cup X \cup \frac{X^2}{2!}
\cup \frac{X^3}{3!} \cup \cdots
\]
at the $X^k/k!=X^k/S_k$ term. The name, however, is our own. There does
not seem to be a satisfactory name in use among geometric 
topologists---indeed, recent authors Mostovoy and Handel do not use
any name at all---and while symmetric product has stayed in use 
among general topologists we prefer to use this for $X^k/S_k$. We therefore
propose the descriptive name ``$k$th finite subset space'' used here and
in our previous paper.

\subsection{Summary of results}

We study the finite subset spaces of a connected graph $\Gamma$ using 
techniques from our previous paper~\cite{circles02} on \exps{k}.
Since $\exp_k$ is a homotopy functor we
may reduce to the case where $\Gamma$ has a single vertex, and accordingly
define $\Gamma_n$ to be the graph with one vertex $v$ and $n$ edges 
$e_1,\ldots,e_n$. Our first result is a complete calculation of the 
homology of \expgb{k}{n}\ and \expg{k}{n}\ for each $k$ and $n$:

\begin{theorem}
\label{Hofexpg.th}
The reduced homology groups of \expgb{k}{n}\ vanish outside dimension
$k-1$ and those of \expg{k}{n} vanish outside dimensions
$k-1$ and $k$. The non-vanishing groups are free. 
The maps 
\[
i\co\expgb{k}{n}\hookrightarrow\expg{k}{n}
\]
and
\[
\cup\{v\}\co\expg{k}{n}\rightarrow\expgb{k+1}{n}
\]
induce isomorphisms on $H_{k-1}$ and $H_k$ respectively while 
\[
\expg{k}{n}\hookrightarrow\expg{k+1}{n}
\]
is twice $(i\circ\cup\{v\})_*$ on $H_k$. The common rank of 
\[
H_k(\expg{k}{n})\cong H_{k}(\expg{k+1}{n})\cong H_{k}(\expgb{k+1}{n})
\]
is given by
\begin{align}
b_k(\expg{k}{n}) &=  \sum_{j=1}^k (-1)^{j-k} \binom{n+j-1}{n-1} \nonumber \\
                 &= \begin{cases}
                    \sum_{j=1}^{\ell}\binom{n+2j-2}{n-2} & 
                                              \text{if $k=2\ell$ is even,} \\
                    n+\sum_{j=1}^{\ell}\binom{n+2j-1}{n-2} &
                                        \text{if $k=2\ell+1$ is odd.}
                    \end{cases}
                    \label{cased_bk.eq}
\end{align}
\end{theorem}

A list of Betti numbers $b_k(\expg{k}{n})$ for $1\leq k\leq 20$ and 
$1\leq n\leq 10$ appears as table~\ref{bettinumbers.tab} in the appendix 
on page~\pageref{bettinumbers.tab}. 

In the case of a circle the homology and fundamental group of \exps{k}\
were enough to determine its homotopy type completely. The argument no longer
applies to \expg{k}{n}, $n\geq 2$, and its applicability for $n=1$ is perhaps 
properly
regarded as being due to a ``small numbers co-incidence'', the vanishing of
$H_{2\ell}(\expg{2\ell}{1})$. However, the argument does apply to
\expgb{k}{n}, and for $k\geq 2$ we have the following:
\begin{theorem}
For $k\geq 2$ the space \expgb{k}{n}\ has the homotopy type of a wedge of 
$b_{k-1}(\expgb{k}{n})$ $(k-1)$--spheres.
\label{wedgeofspheres.th}
\end{theorem}

Having calculated the homology of \expg{k}{n}\ we turn our attention to
the maps $(\Exp{k}\phi)_*$ induced by maps
$\phi\co\Gamma_n\rightarrow\Gamma_m$. Our main result is to reduce the 
problem of calculating such maps to one of finding images of chains under
maps
\begin{align*}
\exps{1} &\rightarrow \exps{1}  \\
\intertext{and}
\exps{2} &\rightarrow \expg{2}{2}
\end{align*}
induced by maps $\s\rightarrow\s$ and $\s\rightarrow\Gamma_2$ respectively.
The reduction is achieved by
defining a ring without unity structure on a subgroup $\tilde\C_*$ of the 
cellular chain complex of $\expinf{\Gamma_n}$. 
The subgroup carries the top homology 
of \expg{k}{n}\ and is preserved by chain 
maps of the form $(\expinf\phi)_\sharp$, and the ring structure is defined 
in such a way that these chain maps are ring homomorphisms. 
The ring $\tilde\C_*\otimes_\integer\rational$ is generated over 
\rational\ by cells in dimensions one and two, leaving a mere
$2n$ cells whose images must be found directly. 

As an application of these results and as an illustration of how much
$(\Exp{k}{\phi})_*$ remembers about $\phi$ we study the action of the
braid group $B_n$ on $H_k(\expg{k}{n})$. The braid group may be
regarded as the mapping class group of a punctured disc and as such it
acts on the graph $\Gamma_n$ via homotopy equivalence. We show that,
for a suitable choice of basis, the braid group acts by block
upper-triangular matrices whose diagonal blocks are representations of
$B_n$ that factor through $S_n$. Consequently, the image of the pure
braid group consists of upper-triangular matrices with ones on the
diagonal and is therefore nilpotent. The number of blocks depends
mildly on $k$ and $n$ but is no more than about $n/2$, and this leads
to a bound on the length of the lower central series.

We remark that the main results of this paper may be used to study the finite 
subset spaces of a closed surface $\Sigma$ via Mayer-Vietoris type arguments.
This may be done by constructing a cover of $\Exp{k}{\Sigma}$ such that
each element of the cover and each $m$--fold intersection is a finite
subset space of a punctured surface, as follows. Choose $k+1$ distinct points
$p_1,\ldots,p_{k+1}$ in $\Sigma$ and let 
$\mathcal{U}_i=\Exp{k}{\bigl(\Sigma\!\setminus\!\{p_i\}\bigr)}$. 
The $\mathcal{U}_i$ form an open cover of $\Exp{k}{\Sigma}$, since each
$\Lambda\in\Exp{k}{\Sigma}$ must omit at least one of the $p_i$, and
moreover each $m$--fold intersection has the form
\[
\bigcap_{j=1}^m \mathcal{U}_{i_j}
      = \Exp{k}{\bigl(\Sigma\!\setminus\!\{p_{i_1},\ldots,p_{i_m}\}\bigr)},
\]
a finite
subset space of a punctured surface as desired. The results of this paper
may then be used to calculate the homology of each intersection and the
maps induced by inclusion, leading to a spectral sequence for
$H_*(\Exp{k}{\Sigma})$.

In~\cite{vanish03} this idea is used to prove two vanishing
theorems for the homotopy and homology groups of the finite subset spaces
of a connected cell complex.

\subsection{Outline of the paper}
\label{outline.sec}

The calculation of the
homology of \expg{k}{n}\ and \expgb{k}{n} is
the main topic of section~\ref{homology.sec}. We find explicit cell
structures for these spaces in section~\ref{cellstructure.sec} and
use them to calculate their fundamental groups. We then show that the reduced
chain complex of \expinf{(\Gamma_n,v)}\ is exact in section~\ref{exact.sec}, 
and use this to prove Theorems~\ref{Hofexpg.th} and~\ref{wedgeofspheres.th} in
section~\ref{Hofexpg.sec}. We give an explicit basis for $H_k(\expg{k}{n})$
in section~\ref{basis.sec} and close with generating functions for 
the Betti numbers $b_k(\expg{k}{n})$ in section~\ref{genfun.sec}. A table
of Betti numbers for $1\leq k\leq 20$ and $1\leq n\leq 10$ appears in
the appendix on page~\pageref{bettinumbers.tab}.

We then turn to the calculation of induced maps in section~\ref{maps.sec}.
The ring structure on $\tilde\C_*$ is motivated and defined in 
section~\ref{chainring.sec} and we show that maps $\phi\co\Gamma_n\rightarrow
\Gamma_m$ induce ring homomorphisms in section~\ref{ringmap.sec}. As 
illustration of the ideas we calculate two examples in 
section~\ref{examples.sec}, the first reproducing a result 
from~\cite{circles02} and the second relating to the generators of the
braid groups. We then state and prove the structure theorem for the
braid group action in sections~\ref{braidintro.sec} and~\ref{nilpotent.sec},
and conclude by looking at the action of $B_3$ on $H_3(\expg{3}{3})$ in
some detail in section~\ref{actiononh3exp3g3.sec}.

\subsection{Notation and terminology}
\label{notation.sec}

We take a moment to fix some language and notation that will be
used throughout.

We will work exclusively with graphs having just one vertex, so 
as above we define $\Gamma_n$ to be the graph with one vertex $v$ and $n$ 
edges $e_1,\ldots,e_n$. Write $I$ for the interval $[0,1]$, 
and for each non-negative integer $m$ let 
$[m]=\{i\in\integer|1\leq i\leq m\}$. We parameterise
$\Gamma_n$ as the quotient of $I\times [n]$
by the subset $\{0,1\}\times [n]$, sending $\{0,1\}\times [n]$ to 
$v$ and $[0,1]\times\{i\}$ to $e_i$. This directs each edge, allowing us
to order any subset of its interior, and we will use this extensively.

Associated to a finite subset $\Lambda$ of $\Gamma_n$ is an $n$--tuple
$\J(\Lambda)=(j_1,\ldots,j_n)$ of non-negative integers
$j_i=|\Lambda \cap \intr e_i|$. Given an $n$--tuple $J=(j_1,\ldots,j_n)$ 
we define its support $\supp J$ to be
\[
\supp J = \{i\in [n]\big| j_i\not=0\}
\] 
and its norm $|J|$ by 
\[
|J|= \sum_{i=1}^n j_i \, .
\]
Note that
\[
|\J(\Lambda)| = \left\{
\begin{array}{ccl}
|\Lambda|   & ~ & \mbox{if $v\not\in\Lambda$}, \\
|\Lambda|-1 & ~ & \mbox{if $v\in\Lambda$}.
\end{array}\right.
\]
In addition we define the mod $2$ support and norm by
\[
\supp_2(J)=\{i\in [n]\big| j_i\not\equiv 0\bmod 2\}
\]
and
\[
|J|_2 = \left|\supp_2(J)\right|.
\]

Bringing two points together in the interior of $e_i$ or moving a point
to $v$ decreases $\J(\Lambda)_i$ by one. It will be convenient to have
some notation for this, so we define 
\[
\drop_i(J) = (j_1,\ldots,j_i-1,\ldots,j_n)
\]
provided $j_i\geq 1$. Lastly, for each subset $S$ of $[n]$ and $n$--tuple
$J$ we write $J|_S$ for the $|S|$--tuple obtained by restricting the index
set to $S$.

\section{The homology of finite subset spaces of graphs}
\label{homology.sec}

\subsection{Introduction}
\label{homologyintro.sec}

Our first step in calculating the homology of a finite subset space
of a connected graph $\Gamma$ is to find explicit cell structures for
\expgb{k}{n}\ and \expg{k}{n}. The approach will be similar to that taken 
in~\cite{circles02}, and we will make use of the boundary map calculated 
there. However, we will adopt a different orientation convention, with the
result that some signs will be changed.

Our cell structure for \expgb{k}{n}\ will consist of 
one $j$--cell $\vcell{J}$ for each $n$--tuple $J$ such that 
$|J|=j\leq k-1$, the interior of $\vcell{J}$ containing those
$\Lambda\in\expgb{k}{n}$ such that $\J(\Lambda)=J$. A cell structure
for \expg{k}{n}\ will be obtained by adding additional cells $\ecell{J}$
for each $J$ with $|J|=j\leq k$; the interior of $\ecell{J}$ will 
contain those $\Lambda\subseteq\Gamma_n\setminus\{v\}$ such that
$\J(\Lambda)=J$. By a ``stars and bars'' argument there are
$\tbinom{n+j-1}{n-1}$ solutions to 
\[
j_1 + \cdots + j_n = j
\]
in non-negative integers (count the arrangements of $j$ ones and $n-1$ 
pluses), so that
\begin{equation}
\chi(\expgb{k}{n}) = \sum_{j=0}^{k-1} (-1)^j \binom{n+j-1}{n-1}
\label{chiofexpgb.eq}
\end{equation}
and
\[
\chi(\expg{k}{n}) =1+2\biggl[\sum_{j=1}^{k-1} (-1)^j \binom{n+j-1}{n-1}\biggr]
                             +(-1)^k\binom{n+k-1}{n-1}.
\]
A cell structure may be found in a similar way for an 
arbitrary connected graph $\Gamma$, with up to $2^{|V(\Gamma)|}$ $j$--cells
for each $|E(\Gamma)|$--tuple $J$ with $|J|=j$.

In these cell structures the spaces \expgb{k+1}{n}\ and \expg{k+1}{n}\ are
obtained from \expgb{k}{n}\ and \expg{k}{n}\ by adding cells in
dimensions $k$ and $k+1$. This has the following consequence for
their homotopy groups. The $(k-1)$--skeleta of
\expgb{k}{n}\ and \expgb{\ell}{n}\ co-incide for $\ell\geq k$, and this
means
that the map on $\pi_i$ induced by the inclusion
$\expgb{k}{n}\hookrightarrow\expgb{\ell}{n}$ is an isomorphism for
$i\leq k-2$. By Handel~\cite{handel00} this map is
zero for $\ell=2k-1$, implying that \expgb{k}{n}\ (and by a
similar argument \expg{k}{n}) is $(k-2)$--connected. It follows 
immediately that the augmented chain complex of $\expinf{(\Gamma_n,v)}$ 
is exact,  and in conjunction with the Euler
characteristic~\eqref{chiofexpgb.eq} and the boundary 
maps~\eqref{bdrysigmaJ.eq} and~\eqref{bdrytildesigmaJ.eq} this is
enough to prove Theorem~\ref{Hofexpg.th}. We nevertheless show directly
that this chain complex is exact in section~\ref{exact.sec} in order
to find bases for the homology groups in section~\ref{basis.sec}.

The fact that the $k$th finite
subset space of a connected graph is $(k-2)$--connected can be used to 
show that
the same conclusion holds for the $k$th finite subset space of a connected
cell complex~\cite{vanish03}. 

\subsection{Cell structures for \expg{k}{n}\ and \expgb{k}{n}}
\label{cellstructure.sec}

We now proceed more concretely. Following the strategy of~\cite{circles02},
each element $\Lambda\in\Exp{j}{e_i}$ has at least one representative
$(x_1,\ldots,x_j)\in[0,1]^j$ such that $x_1\leq \cdots \leq x_j$. Define
simplices
\[
\Delta_j   =    \{(x_1,\ldots,x_{j+1}) | 
                       0\leq x_1 \leq \cdots \leq x_j \leq x_{j+1} = 1 \}
\]
for each $j\geq 0$, and
\[
\tilde\Delta_j           =   \{(x_1,\ldots,x_j) | 
                             0\leq x_1 \leq \cdots \leq x_j \leq 1 \}
\]
for each $j\geq 1$. There are surjections 
$\Delta_j\rightarrow\Exp{j+1}{(e_i,v)}$, 
$\tilde\Delta_j\rightarrow\Exp{j}{e_i}$,
and we let $\vcell[i]j$ be the composition
\[
\Delta_j\rightarrow\Exp{j+1}{(e_i,v)}\hookrightarrow\expgb{j+1}{n},
\]
$\ecell[i]{j}$ the composition
\[
\tilde\Delta_j\rightarrow\Exp{j}{e_i}\hookrightarrow\expg{j}{n}.
\]
We give $\Delta_j$ and $\tilde\Delta_j$ each the orientation
$[x_1,\ldots,x_j]$, a convention that disagrees with the one used
in~\cite{circles02} for some $j$. There $\tilde\Delta_j$ was oriented by
letting its $i$th vertex be 
\[
v_i = (\,\underbrace{0,\ldots,0}_{j-i},\underbrace{1,\ldots,1}_i\,)
\]
for $i=0,\ldots,j$, and the sign of this orientation relative
to the standard one on $\real^j$ is given by
\[
\det [(v_1-v_0)^T,\ldots,(v_j-v_0)^T] =
\begin{cases}
+1 & j\equiv 0,1 \bmod 4, \\
-1 & j\equiv 2,3 \bmod 4.
\end{cases}
\]
A similar calculation shows the same conclusion holds for
$\Delta_j$. To account for this difference we should insert a minus
sign in the boundary map calculated in~\cite{circles02} precisely when
it is applied to $\sigma=\vcell[i]{j}$ or \ecell[i]{j}\ for $j$ even,
since then exactly one of $\sigma$ and $\partial\sigma$ has been given
the opposite orientation. Note however that
$\partial\ecell[i]{j}=\partial\vcell[i]{j}=0$ for $j$ odd, simplifying the
matter and allowing us to simply insert a minus everywhere.

Returning to the discussion at hand, given an $n$--tuple $J$ let
\begin{gather*}
\Delta_J   =  \Delta_{j_1}\times\cdots\times\Delta_{j_n}, \\
\tilde\Delta_J     =  \tilde\Delta_{j_1}\times\cdots\times\tilde\Delta_{j_n},
\end{gather*}
omitting any empty factor $\tilde\Delta_0$ from this last product.
Finally let
\begin{gather*}
\vcell{J} = \vcell[1]{j_1} \cup \cdots \cup \vcell[n]{j_n} \co 
   \Delta_J     \rightarrow \expgb{|J|+1}{n}, \\
\ecell{J} = \ecell[1]{j_1}\cup\cdots\cup\ecell[n]{j_n} \co
   \tilde\Delta_J \rightarrow \expg{|J|}{n},
\end{gather*}
again omitting any factor with $j_i=0$ from $\ecell{J}$.
Each of $\vcell{J}\big|_{\intr\Delta_J}$, 
$\ecell{J}\big|_{\intr\tilde\Delta_J}$ is a homeomorphism of an
open $|J|$--ball onto its image, and we claim:

\begin{lemma}
The spaces \expgb{k}{n}\ and \expg{k}{n}\ have cell structures consisting
respectively of $\{\vcell{J}\big| |J|\leq k-1\}$ and of 
$\{\vcell{J}\big| |J|\leq k-1\}\cup \{\ecell{J}\big| 1\leq|J|\leq k\}$. 
The boundary maps are given by
\begin{equation}
\partial \vcell{J} = -\sum_{i\in\supp J} \frac{1+(-1)^{j_i}}{2} 
         (-1)^{\left| J|_{[i-1]}\right|}\vcell{\drop_i(J)}
\label{bdrysigmaJ.eq}
\end{equation}
and
\begin{equation}
\partial \ecell{J} = \sum_{i\in\supp J} \frac{1+(-1)^{j_i}}{2} 
         (-1)^{\left| J|_{[i-1]}\right|}
         \left(\ecell{\drop_i(J)}-2\vcell{\drop_i(J)}\right).
\label{bdrytildesigmaJ.eq}
\end{equation}
\label{expgcellstructure.lem}
\end{lemma}

Notice that the behaviour of a cell under the boundary map depends only
on the support and parity pattern of $J$. This fact will be of importance
in understanding the chain complexes in section~\ref{exact.sec}.

\begin{proof}
Each element $\Lambda\in\expg{k}{n}$ lies in the interior of the image
of precisely one cell $\vcell{J}$ or $\ecell{J}$, namely
$\vcell{\J(\Lambda)}$ if $v\in\Lambda$ and $\ecell{\J(\Lambda)}$ 
if $v\not\in\Lambda$. The image of $\vcell{J}$ is contained in
\expgb{|J|+1}{n}\ and that of $\ecell{J}$ in \expg{|J|}{n},
so we may set the $j$--skeleton of \expgb{k}{n}\ equal
to \expgb{j+1}{n}\ and the $j$--skeleton of \expg{k}{n}\ equal to 
$(\expg{j}{n})\cup(\expgb{j+1}{n})$ for $j< k$ and \expg{k}{n} for $j=k$.
The boundary of $\tilde\Delta_j$ is found by replacing one or more inequalities
in $0\leq x_1 \leq \cdots \leq x_j \leq 1$ with equalities, resulting in 
fewer points in the interior of $e_i$; thus the image of the 
boundary of $\tilde\Delta_J$ under
$\ecell{J}$ is contained in $\expg{|J|-1}{n}\cup\expgb{|J|}{n}$. 
Similarly, $\vcell{J}$ maps the boundary of $\Delta_J$ into
\expgb{|J|}{n}. So the boundary of a $j$--cell is contained in the 
$(j-1)$--skeleton, and the $\vcell{J}$, $\ecell{J}$ form cell 
structures as claimed.

To calculate the boundary map we use Lemma~1 of~\cite{circles02}, which
with our present notation and orientation convention says
\begin{gather}
\partial\vcell[i]j = -\frac{1+(-1)^j}{2}\vcell[i]{j-1}, \nonumber \\
\partial\ecell[i]{j}=
          \frac{1+(-1)^j}{2}(\ecell[i]{j-1}-2\vcell[i]{j-1}),
\label{bdrytildesigma.eq}
\end{gather}
together with the relation $\partial(\sigma\times \tau)=
(\partial\sigma)\times\tau+(-1)^{\dim\sigma}\sigma\times(\partial\tau)$. 
Calculating the boundary of 
$\ecell[1]{j_1}\times\cdots\times\ecell[1]{j_n}$ and then applying
$\cup_*$ it follows that
\[
\partial \ecell{J} = \sum_{i\in\supp J} (-1)^{\left| J|_{[i-1]}\right|}
 \ecell[1]{j_1}\cup\cdots\cup\partial\ecell[i]{j_i}\cup\cdots
 \cup\ecell[n]{j_n}\, .
\]
Substituting~(\ref{bdrytildesigma.eq}) and observing that 
$\ecell[1]{j_1}\cup\cdots\cup\vcell[i]{j_1}\cup\cdots\cup
\ecell[n]{j_n}=\vcell{J}$ gives~(\ref{bdrytildesigmaJ.eq}), 
and~(\ref{bdrysigmaJ.eq}) follows by a similar argument
or by using $\vcell{J}=(\cup\{v\})_\sharp\ecell{J}$. 
\end{proof}

Let $\C_*$ be the free abelian group generated 
by the $\vcell{J}$, $0\leq |J| < \infty$, and $\tilde\C_*$ the free abelian 
group generated by the $\ecell{J}$, $1\leq |J| < \infty$, each graded by 
degree. Then
\[
H_*(\expgb{k}{n}) = H_*(\C_{\leq k-1})
\]
and
\[
H_*(\expg{k}{n}) = H_*(\C_{\leq k-1}\oplus\tilde\C_{\leq k}).
\]
As discussed at the end of section~\ref{homologyintro.sec} we know
a priori that $(\C_*,\partial)$ is exact except at $\C_0$. We nevertheless
give a direct proof of this in section~\ref{exact.sec}, with a view to 
constructing explicit bases for the homology groups in 
section~\ref{basis.sec} after calculating their ranks in 
section~\ref{Hofexpg.sec}.
Before doing so however we use Lemma~\ref{expgcellstructure.lem} to calculate 
the fundamental groups of \expg{k}{n}\ and \expgb{k}{n} for each $k$ and $n$,
showing directly that 
\expg{k}{n}\ and \expgb{k}{n}\ are simply connected for $k\geq 3$.

\begin{theorem}
The fundamental group of \expg{k}{n}\ is
\begin{enumerate}
\item
free of rank $n$ if $k=1$;
\item
free abelian of rank $n$, containing $i_*\pi_1(\expg{1}{n})$ as a subgroup
of index $2^n$, if $k=2$; and
\item
trivial if $k\geq 3$.
\end{enumerate}
The fundamental group of \expgb{k}{n}\ is free of rank $n$ if $k=2$ and
trivial otherwise.
\end{theorem}

\begin{figure}[tb]
\begin{center}
\leavevmode
\small
\psfrag{tsigi}{$\vcell[i]1$}
\psfrag{tsigj}{$\vcell[j]1$}
\psfrag{sigi}{$\ecell[i]{1}$}
\psfrag{sigij}{$\ecell[i]{1}\cup\ecell[j]{1}$}
\psfrag{sig2i}{$\ecell[i]{2}$}
\psfrag{tsig2i}{$\vcell[i]2$}
\psfrag{a}{(a)}
\psfrag{b}{(b)}
\psfrag{c}{(c)}
\includegraphics{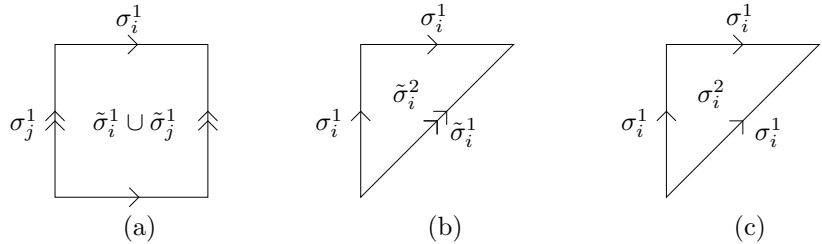}
\caption[Relations in $\pi_1(\expg{k}{n})$ arising from the $2$--cells]
{Relations in $\pi_1(\expg{k}{n})$ arising from the $2$--cells. The
boundary of a cell is found by moving a point in the interior of an 
edge to $v$, or bringing two points in the interior into co-incidence. The
first gives an untilded cell and the second a cell of the same 
kind as the interior. In (a) we
see a torus killing the commutator of $[\vcell[i]1]$ and 
$[\vcell[j]{1}]$; in (b) a M\"obius strip with fundamental group
generated by $[\vcell[i]1]$ and boundary $[\ecell[i]{1}]$; and in
(c) a dunce cap killing $[\vcell[i]1]$.}
\label{relations.fig}
\end{center}
\end{figure}

\begin{proof}
In the unbased case \expg{k}{n}\ the result is obvious for $k=1$ so 
consider $k=2$. The
group $\pi_1(\expg{2}{n})$ is generated by $[\vcell[i]1]$, 
$[\ecell[i]{1}]$, $1\leq i\leq n$, with relations arising from 
the $\ecell[i]{2}$ and $\ecell[i]{1}\cup\ecell[j]{1}$, 
$i\not=j$. The
image of $\ecell[i]{1}\cup\ecell[j]{1}$ is a torus with meridian 
$[\vcell[i]1]$ and longitude $[\vcell[j]{1}]$, while the 
image of $\ecell[i]{2}$ is a M\"obius strip that imposes  the relation
$[\ecell[i]{1}]=[\vcell[i]1]^2$ (see figures~\ref{relations.fig}(a)
and (b)).
It follows that $\pi_1(\expg{2}{n})$ is free abelian with generators
$[\vcell[i]1]$, $1\leq i\leq n$, and that 
$i_*\pi_1(\expg{1}{n})=
\langle [\ecell[1]{1}],\ldots,[\ecell[n]{1}]\rangle$
has index $2^n$. When $k\geq 3$ there are no new generators and additional 
relations $[\vcell[i]1]=1$ from 
each $\vcell[i]2$ (see figure~\ref{relations.fig}(c)), so that
$\pi_1(\expg{k}{n})=\{1\}$.

In the based case $\expgb{1}{n}=\{\{v\}\}$, the map
$\cup\{v\}\co\expg{1}{n}\rightarrow\expgb{2}{n}$ is a homeomorphism, and
for $k\geq 3$ the relations $[\vcell[i]1]=1$ 
from the $\vcell[i]2$ apply as above.
\end{proof}

\subsection{Direct proof of the exactness of $\C_{\geq 1}$}
\label{exact.sec}

We show directly that $\C_{*}$ is exact at each $\ell > 0$ by 
expressing it as a sum of finite subcomplexes and showing that each summand 
is exact. This decomposition will be used to construct explicit bases for
the homology in section~\ref{basis.sec}.

As a first reduction, for each subset $S$ of $[n]$ let $\C^S_*$ be the
free abelian group generated by $\{\vcell{J}|\supp J =S\}$. Since
$\partial\vcell{J}$ is a linear combination of 
the cells $\vcell{\drop_i(J)}$ with $i\in\supp J$ and $
j_i\equiv 0\bmod 2$, each 
$\C^S_*$ is a subcomplex and we have
\[
\C_* = \bigoplus_{S\subseteq [n]} \C^S_* \, .
\]
Note that $\C^\emptyset_* = \C_0$.
Clearly $\C^S_*\cong\C^T_*$ if $|S|=|T|$ so we will show 
$\C^{[m]}_*$ is exact for each $m>0$.

\begin{figure}[tb]
\begin{center}
\small
\leavevmode
\psfrag{(1,2,3)}{$\{1,2,3\}$}
\psfrag{(1,2)}{$\{1,2\}$}
\psfrag{(1,3)}{$\{1,3\}$}
\psfrag{(2,3)}{$\{2,3\}$}
\psfrag{(1)}{$\{1\}$}
\psfrag{(2)}{$\{2\}$}
\psfrag{(3)}{$\{3\}$}
\psfrag{e}{$\emptyset$}
\includegraphics{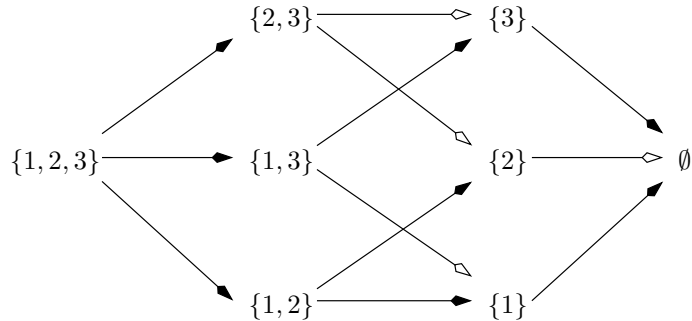}
\caption[The $3$--cube complex $\cu^3_*$]
{The $3$--cube complex $\cu^3_*$. The lattice of subsets of 
$\{1,2,3\}$ forms a $3$--dimensional cube and $\partial S$ is a signed sum of
the neighbours of $S$ of smaller cardinality. In the diagram positive
terms are indicated by solid arrowheads, negative terms by empty 
arrowheads.}
\label{cube.fig}
\end{center}
\end{figure}

We claim that $\C^{[m]}_*$ may be regarded as a sum of many 
copies of a single finite complex, the $m$--cube complex. For each 
$m$--tuple $L$ with all entries odd let $\C^{[m]}_*(L)$ be
the subgroup of $\C^{[m]}_*$ generated by 
$\left\{\vcell{J}\big|j_i-\ell_i\in\{0,1\}\right\}$. Again the 
fact that $\partial\vcell{J}$ is a linear combination of 
$\{\vcell{\drop_i(J)}|i\in\supp J, j_i\equiv 0\bmod 2\}$
implies $\C^{[m]}_*(L)$ is a subcomplex, and moreover that
\[
\C^{[m]}_* = \bigoplus_{L:|L|_2=m} \C^{[m]}_*(L) .
\]
Further, on translating each $m$--tuple by $(\ell_1-1,\ldots,\ell_m-1)$ 
each $\C^{[m]}_*(L)$ can be seen to be isomorphic to 
$\C^{[m]}_*((1,\ldots,1))$ with its grading shifted by $|L|-m$.
We call this common isomorphism class of complex the $m$--cube 
complex $\cu^m_*$, and, replacing $J$ with 
the set of indices of its even entries, will take the free abelian group 
generated by 
the power set of $[m]$, graded by cardinality and with boundary map
\[
\partial S = \sum_{i\in S} (-1)^{\left|[i-1]\setminus S\right|}S\setminus\{i\}
\]
to be its canonical representative; 
for aesthetic purposes we are dropping the minus sign outside the sum. 
The name $m$-cube complex comes from the fact that
the lattice of subsets of $[m]$ forms an $m$--dimensional cube, and that 
$\partial S$ is a signed sum of the neighbours of $S$ of smaller degree. See 
figure~\ref{cube.fig} for the case $m=3$.

Let 
\[
\V_j = \{ S\subseteq [m] \big| |s|=j\mbox{ and } 1\in S\}.
\]
The exactness of $\C^{[m]}_*$  follows from the first statement
of the following lemma; the second statement will be used in 
section~\ref{basis.sec} to construct explicit bases for $H_*(\expg{k}{n})$.

\begin{lemma}
\label{cube.lem}
The $m$--cube complex $\cu^m_*$ is exact. The homology of the truncated
complex $\cu^m_{\leq j}$ is free of rank $\tbinom{m-1}{j}$ in dimension
$j$, with basis $\{\partial S|S\in \V_{j+1}\}$, and zero otherwise.
\end{lemma}

\begin{proof}
We claim that $\V_j\cup \partial \V_{j+1}$ forms a basis for $\cu^m_j$,
from which the lemma follows. Since $\V_j\cup \partial \V_{j+1}$ has at most 
$\tbinom{m-1}{j-1}+\tbinom{m-1}{j}=\tbinom{m}{j}=\operatorname{rank}\cu^m_j$ 
elements we simply check that $\V_j\cup \partial \V_{j+1}$ spans $\cu^m_j$.
It suffices to show that 
$S\in \Span \V_j\cup \partial \V_{j+1}$ for each subset $S$ of size
$j$ not containing $1$; this follows from
\[
\partial (S\cup\{1\}) = S - 
    \sum_{i\in S} (-1)^{|[i-1]\setminus S|}S\cup\{1\}\setminus \{i\}
\]
if $1\not\in S$. 
\end{proof}

\subsection{The homology groups of \expgb{k}{n}\ and \expg{k}{n}}
\label{Hofexpg.sec}

We calculate the homology groups of \expgb{k}{n}\ and \expg{k}{n}\ 
using the exactness of $\C_{\geq 1}$, the Euler 
characteristic~\eqref{chiofexpgb.eq}, and the 
boundary maps~\eqref{bdrysigmaJ.eq} 
and~\eqref{bdrytildesigmaJ.eq}. Explicit bases are found in 
section~\ref{basis.sec} using the decomposition of $\C_*$ into 
subcomplexes. 

\begin{proof}[Proof of Theorem~\ref{Hofexpg.th}]
Since \expgb{k}{n}\ is path connected with homology equal to that of
$\C_{\leq {k-1}}$, its reduced homology 
vanishes except perhaps in 
dimension $k-1$, by the exactness of $\C_{\geq 1}$. Moreover $H_{k-1}$ is 
equal to $\ker \partial_{k-1}$ and is therefore free; its rank may be found
using $\chi(\expgb{k}{n}) = b_0 + (-1)^{k-1}b_{k-1}$
and equation~\eqref{chiofexpgb.eq}, yielding
\[
b_{k-1}(\expgb{k}{n}) = \sum_{j=1}^{k-1} (-1)^{k-j-1} \binom{n+j-1}{n-1}.
\]
This may be expressed as a sum of purely positive terms by grouping the
summands in pairs, starting with the largest, and using 
$\tbinom{p}{q}-\tbinom{p-1}{q}=\tbinom{p-1}{q-1}$. Doing this for
$b_k(\expgb{k+1}{n})$ gives the expression in equation~\eqref{cased_bk.eq}.

Now consider $H_*(\expg{k}{n})=H_*(\C_{\leq{k-1}}\oplus\C_{\leq{k}})$.
Write  $C_j$ for the $j$th chain group of 
$\C_{\leq{k-1}}\oplus\tilde\C_{\leq{k}}$, $\Z_j$ for the $j$--cycles in 
$\C_j$ and $Z_j$ for the $j$--cycles in $C_j$.
Extending $\vcell{J}\mapsto\ecell{J}$ linearly to a group isomorphism from
$\C_j$ to $\tilde\C_j$ for each $j\geq 1$, the boundary 
maps~\eqref{bdrysigmaJ.eq} and~\eqref{bdrytildesigmaJ.eq} give
\[
\partial\tilde c = 2\partial c - \widetilde{\partial c}
\]
for each chain $c\in \C_{\geq 1}$. It
follows that $Z_j=\Z_j\oplus\tilde\Z_j$ for $1\leq j\leq k-1$ and 
that $H_k(\expg{k}{n})=Z_k$ is equal to $\tilde\Z_k$. Moreover 
\[
\partial C_k = \{2z-\tilde{z} | z\in\Z_{k-1}\}
\]
by the exactness of $\C_*$ at $\C_{k-1}$, 
so that $H_{k-1}(\expg{k}{n})=\Z_{k-1}$. We show that the remaining 
reduced homology groups vanish.

Let $z=z_1\oplus\tilde z_2\in Z_j$ for some $1\leq j\leq k-2$.
By the exactness of 
$\C_{\geq 1}$ there are $w_1,w_2\in\C_{j+1}$ such that 
\begin{align*}
\partial w_1 & = z_1 +2z_2, \\
\partial w_2 & = z_2.
\end{align*}
Since $j\leq k-2$ we have $w_1-\tilde{w}_2\in C_{j+1}$, and
\begin{align*}
\partial(w_1-\tilde{w}_2) 
                  &= \partial w_1 - \partial \tilde{w}_2 \\
                  &= z_1 + 2z_2 - 2\partial w_2 + \widetilde{\partial w_2} \\
                  &= z_1 + \tilde{z}_2,
\end{align*}
so that $C_*$ is exact at $C_j$ as claimed.

It remains to determine the maps induced by
\begin{gather*}
i\co\expgb{k}{n}\hookrightarrow\expg{k}{n}, \\
\cup\{v\}\co\expg{k}{n}\rightarrow\expgb{k+1}{n} 
\end{gather*}
and
\[
\expg{k}{n}\hookrightarrow\expg{k+1}{n}
\]
on homology. In each case there is only one dimension in which the induced map 
is not trivially zero. We have 
$H_{k-1}(\expgb{k}{n})=\Z_{k-1}=H_{k-1}(\expg{k}{n})$, so that $i_*$ is an
isomorphism on $H_{k-1}$, and adding $v$ to each element of \expg{k}{n}\ sends
$\tilde z$ to $z$ for each $z\in\Z_k$, inducing an isomorphism on
$H_k$. Lastly, $\expg{k}{n}\hookrightarrow\expg{k+1}{n}$ sends
$[\tilde{z}]$ to $[\tilde{z}]=2[z]=2(i\circ\cup\{v\})_*[\tilde{z}]$ for each 
$z\in\Z_k$, inducing two times $(i\circ\cup\{v\})_*$ as claimed.
\end{proof}

The homology and fundamental group of \expgb{k}{n}\ are enough to determine
its homotopy type completely. When $k=1$ it is a single point $\{\{v\}\}$,
and when $k\geq 2$ we have:

\begin{corollary}[Theorem~\ref{wedgeofspheres.th}]
For $k\geq 2$ the space \expgb{k}{n}\ has the homotopy type of a
wedge of $b_{k-1}(\expgb{k}{n})$ $(k-1)$--spheres.
\end{corollary}
\begin{proof}
Since \expgb{2}{n}\ is homeomorphic to $\Gamma_n$ we may assume $k\geq 3$. But
then \expgb{k}{n}\ is a simply connected Moore space 
$M(\integer^{b_{k-1}},k-1)$
and the result follows from the Hurewicz and Whitehead theorems.
\end{proof}

\subsection{A basis for $H_k(\expg{k}{n})$}
\label{basis.sec}

We use the decomposition of $\C_*$ as a sum of subcomplexes to give 
an explicit basis for $H_k(\expg{k}{n})$.

\begin{theorem}
The set
\[
\mathcal{B}(k,n) = \left\{\widetilde{\partial\vcell{J}} \big| 
  |J|=k+1\mbox{ and } j_i\equiv 0\bmod{2}
  \mbox{ for }i=\min(\supp J)\right\}
\]
is a basis for $H_k(\expg{k}{n})$.
\label{basis.th}
\end{theorem}

\begin{proof}
It suffices to find a basis for $\Z_k$ and map it across to $\tilde\Z_k$.
Extending notation in obvious ways we have
\begin{align*}
\Z_k &= \bigoplus_{S\subseteq[n]} \Z^{S}_k \\
     &= \bigoplus_{S\subseteq[n]} \,
         \bigoplus_{\substack{L:\supp L=S,\\ |L|_2=|S|}}\Z^{S}_k(L).
\end{align*}
Each $\Z^{S}_k(L)$ in this sum is isomorphic to $\mathsf{Z}^{|S|}_j$ for some
$j$, and tracing back through this isomorphism we see that $\V_{j+1}$ is
carried up to sign to 
\[
\mathcal{V}^S_{k+1}(L)=\{\vcell{J}\big| 
|J|=k+1,j_i-\ell_i\in\{0,1\}, j_i\equiv0\bmod 2 \mbox{ for } i=\min(\supp J)\}.
\]
By Lemma~\ref{cube.lem} 
$\{\partial\vcell{}|\vcell{}\in\mathcal{V}^S_{k+1}(L)\}$
is a basis for $\Z^{S}_k(L)$, and taking the union over $S$ and $L$
completes the proof.
\end{proof}

As an exercise in counting we check that $\mathcal{B}(k,n)$ has the right
cardinality. This is equivalent to showing that the number $s(k,n)$
of non-negative integer solutions to 
\[
j_1 + \cdots + j_n = k
\]
in which the first non-zero summand is odd is given by 
equation~\eqref{cased_bk.eq}.
We do this by induction on $k$, inducting separately over the even and
odd integers.

In the base cases $k=1$ and~$2$ there are clearly $n$ and $\tbinom{n}{2}$
solutions respectively. It therefore suffices to show that
$s(k,n)-s(k-2,n)=\tbinom{n+k-2}{n-2}$. Adding two to the first non-zero 
summand gives an injection from solutions  with $k=\ell-2$ to solutions
with $k=\ell$, hitting all solutions except those for which the first
non-zero summand is one. If $j_{n-i}=1$ is the first non-zero summand
then what is left is an unconstrained non-negative integer solution to
\[
j_{n-i+1}+\cdots+j_n = \ell -1,
\]
of which there are $\binom{\ell+i-2}{\ell-1}$, so that
\[
s(k,n)-s(k-2,n)=\sum_{i=1}^{n-1} \binom{k+i-2}{k-1}.
\]
This is a sum down a diagonal of Pascal's triangle and as such is easily
seen to equal $\tbinom{k+n-2}{k}=\tbinom{k+n-2}{n-2}$.

\subsection{Generating functions for $b_k(\expg{k}{n})$}
\label{genfun.sec}

We conclude this section by giving generating functions for the
Betti numbers $b_k(\expg{k}{n})$.

\begin{theorem}
The Betti number $b_k(\expg{k}{n})$ is the co-efficient of $x^k$ in
\begin{equation}
\frac{1-(1-x)^n}{(1+x)(1-x)^n}.
\label{genfun.eq}
\end{equation}
\label{genfun.th}
\end{theorem}

\begin{proof}
The co-efficient of $x^j$ in 
\[
\frac{1}{(1-x)^n} = (1+x+x^2+x^3+\cdots)^n
\]
is the number of non-negative solutions to $j_1+\cdots+j_n=j$, in other
words $\binom{n+j-1}{n-1}$. Multiplication by
$(1+x)^{-1} = 1-x+x^2-x^3+\cdots$
has the effect of taking alternating sums of co-efficients, so we
subtract $1$ first to remove the unwanted degree zero 
term from $(1-x)^{-n}$, arriving at~\eqref{genfun.eq}.
\end{proof}

\section{The calculation of induced maps}
\label{maps.sec}

\subsection{Introduction}
\label{mapsintro.sec}

Given a pointed map $\phi\co(\Gamma_n,v)\rightarrow(\Gamma_m,v)$ there
are induced maps 
\[
(\Exp{k}{\phi})_* \co H_{k-1}(\expgb{k}{n})\rightarrow H_{k-1}(\expgb{k}{m})
\]
and
\[
(\Exp{k}{\phi})_* \co H_{p}(\expg{k}{n})\rightarrow H_{p}(\expg{k}{m})
\]
for $p=k-1,k$. In view of the commutative diagrams
\[
\begin{CD}
\expgb{k}{n} @>{\Exp{k}{\phi}}>> \expgb{k}{m} \\
@VViV                             @VViV       \\
\expg{k}{n}  @>{\Exp{k}{\phi}}>> \expg{k}{m}
\end{CD}
\]
and
\[
\begin{CD}
\expg{k-1}{n} @>{\Exp{k-1}{\phi}}>> \expg{k-1}{m} \\
@VV{\cup\{v\}}V                    @VV{\cup\{v\}}V   \\
\expgb{k}{n}  @>{\Exp{k}{\phi}}>> \expgb{k}{m}
\end{CD}
\]
and the isomorphisms induced by $i$ and $\cup\{v\}$ on $H_{k-1}$ and 
$H_k$ respectively it suffices to understand just one of these, and we 
will focus our attention on $H_k(\expg{k}{n})\rightarrow H_k(\expg{k}{m})$.
The purpose of this section is to reduce the problem of calculating this
map to the problem of finding the images of the basic cells 
$\ecell[i]{1},\ecell[i]{2}$
under the chain map $(\expinf{\phi})_\sharp$. The reduction will be done by
defining a multiplication on $\tilde\C_*$, giving it the structure of
a ring without unity generated by the $\ecell[i]{j}$. The multiplication
will be defined in such a way that the cellular chain map 
$(\expinf{\phi})_\sharp|_{\tilde\C_*}$
is a ring homomorphism, reducing calculating $(\Exp{k}{\phi})_*$ to 
calculations in
the chain ring once the $(\expinf{\phi})_\sharp\ecell[i]{j}$ are found.
The reduction to just the cells $\ecell[i]{1},\ecell[i]{2}$ is achieved
by working over the rationals, as $\tilde\C_*\otimes_\integer\rational$
will be generated over \rational\ by just these $2n$ cells.

In what follows we will assume that $\phi$ is smooth, in the sense
that $\phi$ is smooth on the open set
$\phi^{-1}(\Gamma_m\setminus\{v\})$.  This ensures that \Exp{j}{\phi}\
is smooth off the preimage of the $(j-1)$--skeleton
$(\expg{j}{n})^{j-1}=(\expg{j-1}{n})\cup(\expgb{j}{n})$, allowing us to
use smooth techniques on the manifold
$\expg{j}{n}\setminus(\expg{j}{n})^{j-1}$. Smoothness of $\phi$ may be ensured 
by homotoping it to a standard form defined as follows. The
restriction of $\phi$ to each edge $e_i$ is an element of
$\pi_1(\Gamma_m,v)$, and as such is equivalent to a reduced word $w_i$
in the $\{e_a\}\cup\{\bar{e}_a\}$. We consider $\phi$ to be in standard
form if $\phi|_{e_i}$ traverses each letter of $w_i$ at constant
speed.

\subsection{The chain ring}
\label{chainring.sec}

We observe that the operation $(g,h)\mapsto g\cup h$ suggests
a natural way of multiplying cells and we study it with an eye to
applying the results to maps of the form $(\Exp{k}{\phi})\circ\ecell{J}$.

A map of pairs $g\co 
(B^j,\partial B^j)\rightarrow \bigl(\expg{j}{n},(\expg{j}{n})^{j-1}\bigr)$ 
induces a map
\[
g_*\co 
H_j(B^j,\partial B^j)\rightarrow H_j(\expg{j}{n},(\expg{j}{n})^{j-1}),
\]
and the homology group on the right is canonically isomorphic to the cellular
chain group $\tilde\C_j$. 
Writing $\epsilon^j$ for the positive generator of $H_j(B^j,\partial B^j)$,
if $g$ is smooth on the
open set $g^{-1}(\expg{j}{n}\setminus (\expg{j}{n})^{j-1})$ then
this map is given by
\[
g_* \epsilon^j = \sum_{|J|=j} \pairing{g}{\ecell{J}}\ecell{J}, 
\]
in which $\pairing{g}{\ecell{J}}$ is the
signed sum of preimages of a generic point in the interior of $\ecell{J}$. 
If $h$ is a second map of pairs $(B^\ell,\partial B^\ell)\rightarrow 
\bigl(\expg{\ell}{n},(\expg{\ell}{n})^{\ell-1}\bigr)$ then $g\cup h$ is a map
of pairs
\[
(B^{j+\ell},\partial B^{j+\ell})\rightarrow 
\bigl(\expg{j+\ell}{n},(\expg{j+\ell}{n})^{j+\ell-1}\bigr)
\]
also. The following lemma shows that 
$\cup$ behaves as might be hoped on the chain level.

\begin{lemma}
Given two maps of pairs 
$g\co(B^j,\partial B^j)\rightarrow (\expg{j}{n},(\expg{j}{n})^{j-1})$, 
$h\co(B^\ell,\partial B^\ell)\rightarrow 
(\expg{\ell}{n},(\expg{\ell}{n})^{\ell-1})$,
each smooth off the preimage of the codimension one skeleton, we have
\begin{equation}
{(g\cup h)}_* = \sum_{|J|=j,|L|=\ell} \pairing{g}{\ecell{J}} 
    \pairing{h}{\ecell{L}} {(\ecell{J}\cup\ecell{L})}_*.
\label{chainovercup.eq}
\end{equation}
\label{chainovercup.lem}
\end{lemma}

\begin{proof}
Fix an $n$--tuple $M$ such that $|M|=j+\ell$ and let $\Lambda$ be a generic
point in the interior of $\ecell{M}$. It suffices to check that $g\times h$
and $\sum_{J,L}\pairing{g}{\ecell{J}}\pairing{h}{\ecell{L}}
\ecell{J}\times\ecell{L}$ have the same signed sum of preimages at each
point $(\Lambda',\Lambda'')\in \expg{|J|}{n}\times\expg{|L|}{n}$ such that
$\Lambda'\cup\Lambda''=\Lambda$; note that for cardinality reasons
$(\Lambda',\Lambda'')$ forms a partition of $\Lambda$. For $g\times h$
this signed sum is given by
\begin{align*}
\pairing{g\times h}{\ecell{\J(\Lambda')}\times\ecell{\J(\Lambda'')}}
   &= \sum_{p\in g^{-1}(\Lambda')}\,\sum_{q\in h^{-1}(\Lambda'')}
            \sign(\det D(g\times h)(p,q)) \\
   &= \sum_{p\in g^{-1}(\Lambda')}\,\sum_{q\in h^{-1}(\Lambda'')}
            \sign\bigl(\det Dg(p)\bigr)\sign\bigl(\det Dh(q)\bigr) \\
   &= \biggl[\sum_{g^{-1}(\Lambda')}\sign(\det Dg(p)\bigr)\biggr]
      \biggl[\sum_{h^{-1}(\Lambda'')}\sign(\det Dh(q)\bigr)\biggr] \\
   &= \pairing{g}{\ecell{\J(\Lambda')}}\pairing{h}{\ecell{\J(\Lambda'')}}.
\end{align*}
The lemma follows from the fact that  $\pairing{\ecell{J}\times\ecell{L}}
{\ecell{\J(\Lambda')}\times\ecell{\J(\Lambda'')}}$ is zero unless
$J=\J(\Lambda')$ and $L=\J(\Lambda'')$, in which case it is one.
\end{proof}

Since ${(\ecell{J}\cup\ecell{L})}_*\epsilon^{j+\ell}$ is a multiple of 
$\ecell{J+L}$
the next step is to understand the pairings 
$\pairing{\ecell{J}\cup\ecell{L}}{\ecell{J+L}}$. 
Interchanging adjacent factors $\ecell[a]{r}$ and $\ecell[b]{s}$ in the 
product $\ecell{J}\cup\ecell{L}$ simply introduces a sign $(-1)^{rs}$, so 
we may gather basic cells from the same edge together and consider 
pairings of the form 
\[
\pairing{(\ecell[1]{j_1}\cup\ecell[1]{\ell_1})\cup\cdots\cup
        (\ecell[n]{j_n}\cup\ecell[n]{\ell_n})}{\ecell{J+L}}
= \prod_{i=1}^{n} \pairing{\ecell[i]{j_i}\cup\ecell[i]{\ell_i}}
                         {\ecell[i]{j_i+\ell_i}}.
\]
The quantity 
$\pairing{\ecell[a]{r}\cup\ecell[a]{s}}{\ecell[a]{r+s}}$
is equal to $\sbinom{r+s}{r}$, 
the $q$--binomial co-efficient
$\qbinom{r+s}{r}$ specialised to $q=-1$. The
correspondence can be seen as follows. Take $r+s$ objects, numbered from 
$1$ to $r+s$ and laid out in order, and paint $r$ of them blue and the 
rest red. Shuffle them so that the blue ones are at the front in ascending
order, followed by the red ones in ascending order, giving an element of
the symmetric group $S_{r+s}$. Then $\sbinom{r+s}{r}$ is the number
of ways of choosing $r$ objects from $r+s$ in this way, counted with the
sign of the associated permutation, and is equal to 
$\pairing{\ecell[a]{r}\cup\ecell[a]{s}}{\ecell[a]{r+s}}$: the blue and 
red points represent the
elements of a generic point in $\Exp{r+s}{e_a}$ coming from $\ecell[a]{r}$
and $\ecell[a]{s}$ respectively, and the derivative at this preimage is
the matrix of the associated permutation.

\begin{figure}
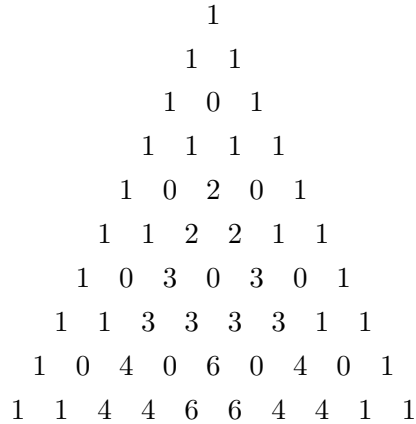

\begin{gather*}
1 \\
1 \quad 1 \\
1\quad 0 \quad 1 \\
1\quad 1\quad 1\quad 1 \\
1\quad 0\quad 2\quad 0\quad 1\\
1\quad 1\quad 2\quad 2\quad 1\quad 1\\
1\quad 0\quad 3\quad 0\quad 3\quad 0\quad 1\\
1\quad 1\quad 3\quad 3\quad 3\quad 3\quad 1\quad 1\\
1\quad 0\quad 4\quad 0\quad 6\quad 0\quad 4\quad 0\quad 1\\
1\quad 1\quad 4\quad 4\quad 6\quad 6\quad 4\quad 4\quad 1\quad 1
\end{gather*}
\caption[Pascal's other triangle]
{The first ten rows of Pascal's other triangle, enough so that the
pattern is clear. Each row of Pascal's triangle appears twice; on the first
occurrence zeros are inserted between each entry, and on the second 
each entry appears twice.}
\label{pascals.fig}
\end{figure}

The calculation of $\sbinom{r+s}{r}$ is the subject of the following lemma.
The result, which we might call ``Pascal's other triangle''---being a much less
popular model than the one we know and love---appears in
figure~\ref{pascals.fig}. For further information on $q$--binomial 
co-efficients and related topics see for example Stanley~\cite{stanley-ec}, 
Kac and Cheung~\cite{kac-cheung-qc}, and
Baez~\cite[weeks 183--188]{baez-twf}.

\begin{lemma}
The value of the signed binomial co-efficient $\sbinom{m}{r}$ is given by
\begin{equation}
\sbinom{m}{r} = 
  \frac{1+(-1)^{r(m-r)}}{2}\binom{\lfloor m/2\rfloor}{\lfloor r/2\rfloor}.
\label{pascals.eq}
\end{equation}
\label{pascals.lem}
\end{lemma}

\begin{proof}
The $q$--binomial co-efficient $\qbinom{m}{r}$ satisfies
$\qbinom{m}{0}=\qbinom{m}{m}=1$ and 
\[
\qbinom{m}{r} = \qbinom{m-1}{r-1} + q^r \qbinom{m-1}{r},
\]
both easily seen for $q=-1$ from the definition above: the recurrence
relation is proved analogously to the familiar one for Pascal's triangle,
the sign $(-1)^r$ arising from shuffling the first object into place if it 
is red instead of blue. The equality~\eqref{pascals.eq} is then readily proved
by induction on $m$, by considering in turn the four possibilities for the 
parities of $m$ and $r$.
\end{proof}

We remark that the equality $\qbinom{r+s}{r}=\qbinom{r+s}{s}$ can be
seen for $q=-1$ from the fact that 
$\lfloor r/2 \rfloor + \lfloor s/2 \rfloor =
\lfloor (r+s)/2 \rfloor$ unless both $r$ and $s$ are odd, in which case
the co-efficient $\sbinom{r+s}{r}$ vanishes. Additionally  
$\sbinom{r+s}{r}\sbinom{r+s+t}{r+s}$
and $\sbinom{r+s+t}{r}\sbinom{s+t}{s}$ are both equal to
\[
\frac{\lfloor (r+s+t)/2\rfloor !}
{\lfloor r/2 \rfloor ! \lfloor s/2 \rfloor ! \lfloor t/2 \rfloor !}
\]
if no more than one of $r,s,t$ is odd, and zero otherwise.

We now define the chain ring of \expinf{\Gamma_n}\ to be the ring without unity
generated over $\integer$ by the set $\{\ecell[i]{j}\big|1\leq i\leq n, 
j\geq 1\}$, with relations 
\begin{gather*}
\ecell[a]j\ecell[b]\ell = (-1)^{j\ell}\ecell[b]\ell\ecell[a]j \\
\ecell[a]j\ecell[a]\ell = {\textstyle\sbinom{j+\ell}{j}}\ecell[a]{j+\ell}
\end{gather*}
for all $1\leq a,b\leq n$ and $j,\ell\geq 1$. 
This definition is again a specialisation to $q=-1$ of a construction that
applies more generally, and is chosen so that the 
conclusion~\eqref{chainovercup.eq} of 
Lemma~\ref{chainovercup.lem} may be rewritten as
\begin{equation}
{(g\cup h)}_*\epsilon^{j+\ell} = (g_*\epsilon^j)(h_*\epsilon^\ell),
\label{chainovercup2.eq}
\end{equation}
in which the multiplication on the right hand side takes place in the
chain ring. The lack of an identity could be easily remedied but we
have chosen not to so that all elements of the ring are chains.
We shall denote the chain ring simply by $\tilde\C_*$, or 
$\tilde\C_*(\Gamma_n)$ in case of ambiguity.

\subsection{Calculating induced maps}
\label{ringmap.sec}

We now have all the machinery required to state and prove 
our main result on the 
calculation of $(\Exp{k}{\phi})_*$, namely that the cellular chain map
$(\expinf{\phi})_\sharp|_{\tilde\C_*}$ is a ring homomorphism. 
This reduces the calculation of the chain map to the calculation
of the images of the basic cells $\ecell[i]{j}$, each an exercise in 
counting points with signs, and multiplication and addition in the
chain ring. 

\begin{theorem}
If $\phi\co(\Gamma_n,v)\rightarrow(\Gamma_m,v)$ is smooth on 
$\phi^{-1}(\Gamma_m\setminus\{v\})$ then the cellular chain map
$(\expinf{\phi})_\sharp\co\tilde\C_*(\Gamma_n)\rightarrow\tilde\C_*(\Gamma_m)$
is a ring homomorphism.
\label{ringmap.th}
\end{theorem}

\begin{proof}
If $\phi$ is smooth off $\phi^{-1}(v)$ then as
noted at the end of section~\ref{mapsintro.sec} the map 
$\Exp{j}{\phi}=\expinf{\phi}|_{\expg{j}{n}}$ is smooth off the codimension
one skeleton for each $j$, so we may use Lemma~\ref{chainovercup.lem}. 
Thus
\begin{align*}
(\expinf{\phi})_\sharp (\ecell{J}\ecell{L}) 
 &= (\Exp{j+\ell}{\phi})_*((\ecell{J}\cup\ecell{L})_*\epsilon^{j+\ell})  && 
     \text{by~\eqref{chainovercup2.eq}}                                 \\
 &= {\bigl((\Exp{j+\ell}{\phi})\circ{(\ecell{J}\cup\ecell{L})}\bigr)}_*
       \epsilon^{j+\ell}                                                 &&  \\
 &= {\bigl(((\Exp{j}{\phi})\circ\ecell{J})
         \cup((\Exp{\ell}{\phi})\circ\ecell{L})\bigr)}_*\epsilon^{j+\ell} && \\
 &= \bigl({\smash{((\Exp{j}{\phi})\circ\ecell{J})}}_*\epsilon^j\bigr)
      \bigl(\smash[b]{((\Exp{\ell}{\phi})\circ\ecell{L})}_*\epsilon^\ell\bigr)
      && \text{by~\eqref{chainovercup2.eq}}                                 \\
 &= \bigl({(\Exp{\smash[b]{j}}{\phi})}_*{\ecell[*]{J}}\epsilon^j\bigr)
        \bigl({(\Exp{\ell}{\phi})}_*{\ecell[*]{L}}\epsilon^\ell\bigr)   &&  \\
 &= \bigl((\expinf\phi)_\sharp\ecell{J}\bigr)
    \bigl((\expinf\phi)_\sharp\ecell{L}\bigr),  &&
\end{align*}
from which the result follows.
\end{proof}

As an immediate consequence we have
\[
(\expinf{\phi})_\sharp \ecell{J} = 
         \prod_{i=1}^n (\expinf{\phi})_\sharp \ecell[i]{j_i},
\]
so that $(\expinf{\phi})_\sharp \ecell{J}$ may be found knowing just
the images of the basic cells \ecell[i]{j}\ as claimed. 
To reduce the number of cells $\ecell{}$ for which
$(\expinf{\phi})_\sharp\ecell{}$ must be calculated directly even further,
observe that 
\begin{equation}
\ecell[i]{j} =
\begin{cases}
\frac{1}{\ell!}\left(\ecell[i]{2}\right)^\ell & \text{if $j=2\ell$}, \\
\frac{1}{\ell!}\ecell[i]{1}\left(\ecell[i]{2}\right)^\ell
        & \text{if $j=2\ell+1$},
\end{cases} 
\label{gensoverQ.eq}
\end{equation}
so that $\tilde\C_*\otimes_\integer\rational$ is generated over \rational\
by the set $\{\ecell[i]{j}\big| 1\leq i\leq n, j=1,2\}$, subject only 
to the relations that the $\ecell[i]{1}$ anti-commute with each other and
the $\ecell[i]{2}$ commute with everything. In essence this
reduces calculating $(\expinf{\phi})_\sharp$ to understanding the 
behaviour of chains under maps
\begin{align*}
\exps{1} &\rightarrow \exps{1}  \\
\intertext{and}
\exps{2} &\rightarrow \expg{2}{2}
\end{align*}
induced by maps $\s\rightarrow\s$ and $\s\rightarrow\Gamma_2$ respectively.
These are both simple exercises in counting points with signs and we give 
the answers, which are easily checked. Let $w$ be a reduced word in 
$\{e_1,e_2,\bar{e}_1,\bar{e}_2\}$, and let $\phi$ from $\s=\Gamma_1$
to $\Gamma_2$ send $e_1$ to $w$. 
Then $\pairing{(\expinf\phi)\circ\ecell[1]{1}}{\ecell[i]{1}}$ and
$\pairing{(\expinf\phi)\circ\ecell[1]{2}}{\ecell[i]{2}}$ are both given by the
winding number of $w$ around $e_i$, and 
$\pairing{(\expinf\phi)\circ\ecell[1]{2}}{\ecell{(1,1)}}$ 
is the number of pairs
of letters $(a_1,a_2)$ in $w$, $a_i\in\{e_i,\bar{e}_i\}$, counted with the
product of a minus one for each bar and a further minus one if $a_2$ occurs 
before $a_1$.

\subsection{Examples}
\label{examples.sec}

As an illustration of the ideas in this section we calculate two examples.
The first reproduces a result from~\cite{circles02} on maps $\s\rightarrow\s$,
and the second will be useful in understanding the action of the braid group.

Let $\phi\co\s\rightarrow\s$ be a degree $d$ map. By Theorem~\ref{Hofexpg.th}
of this paper or Theorem~4 of~\cite{circles02} we have
\[
H_k(\exps{k}) \cong  \begin{cases}
                     0        & \text{$k$ even}, \\
                     \integer & \text{$k$ odd},
                     \end{cases}
\]
so the only map of interest is 
$(\Exp{2\ell-1}{\phi})_*$ on $H_{2\ell-1}$. The homology group 
$H_{2\ell-1}(\exps{2\ell-1})$ is generated by $\ecell[1]{2\ell-1}$, and 
by~\eqref{gensoverQ.eq} and the discussion at the end of 
section~\ref{ringmap.sec} we have
\begin{align*}
(\Exp{2\ell-1}{\phi})_\sharp\ecell[1]{2\ell-1}
  &= \frac{1}{(\ell-1)!}(\Exp{2\ell-1}{\phi})_\sharp
        \bigl(\ecell[1]{1}\left(\ecell[1]{2}\right)^{\ell-1}\bigr) \\
  &= \frac{1}{(\ell-1)!}(d\ecell[1]{1})(d\ecell[1]{2})^{\ell-1}   \\
  &= \frac{d^\ell}{(\ell-1)!}\ecell[1]{1}\left(\ecell[1]{2}\right)^{\ell-1} \\
  &=  d^\ell \ecell[1]{2\ell-1}.
\end{align*}
Thus \Exp{2\ell-1}\phi\ is a degree $d^\ell$ map, as found in Theorem~7
of~\cite{circles02}.

For the second example consider the map $\tau\co\Gamma_2\rightarrow\Gamma_2$
sending $e_1$ to $e_2$ and $e_2$ to $e_2 e_1 \bar{e}_2$. We shall
compare this with the map $\perm\tau$ that simply switches $e_1$ and $e_2$.
These two maps are homotopic through a homotopy that drags
$v$ around $e_2$, so we expect the same induced map once we pass
to homology, but an understanding of the chain maps will be useful in
section~\ref{braid.sec} when we study the action of the braid group.
Clearly $(\expinf{\perm{\tau}})_\sharp$ simply interchanges $\ecell[1]{j}$ and
$\ecell[2]{j}$, and likewise $(\expinf\tau)_\sharp\ecell[1]{j}=\ecell[2]{j}$ 
for each $j$. The only difficulty therefore is in finding 
$(\expinf\tau)_\sharp\ecell[2]{j}$, and we shall
do this in two ways, first by calculating it directly and then by
working over \rational\ using~\eqref{gensoverQ.eq}.

To calculate $(\expinf\tau)_\sharp\ecell[2]{j}$ directly let 
$p=(x_1,\ldots,x_\ell,y_1,\ldots,y_m)$ be a generic point in 
$\Delta_\ell\times\Delta_m$, $\ell+m=j$, $m\geq 1$. Preimages of $p$ come in 
pairs that differ only in whether the preimage of $y_i$ comes from
the letter $e_2$ or $\bar{e}_2$ of $w$. In the first case the 
$(\ell+1)$th row of the derivative has a $3$ in the first column and
zeros elsewhere, and in the second it has a $-3$ in the last column
and zeros elsewhere. Thus each pair contributes $(-1)^\ell+(-1)^m$ times
the sign of the corresponding preimage in $\ecell[2]{j-1}$ of 
$(x_1,\ldots,x_\ell,y_2,\ldots,y_m)\in\Delta_\ell\times\Delta_{m-1}$,
and consequently
\[
\pairing{(\expinf{\tau})\circ\ecell[2]{j}}{\ecell{(\ell,m)}}
     =\bigl((-1)^\ell + (-1)^m\bigr)
       \pairing{(\expinf{\tau})\circ\ecell[2]{j-1}}{\ecell{(\ell,m-1)}}.
\]
This recurrence relation is easily solved to give 
\[
\pairing{(\expinf{\tau})\circ\ecell[2]{j}}{\ecell{(\ell,m)}}=
\begin{cases}
 1  &  \text{if $m=0$,} \\
-2  &  \text{if $\ell$ odd, $m=1$}, \\
 0  &  \text{otherwise},
\end{cases}
\]
so that
\begin{equation}
(\expinf\tau)_\sharp\ecell[2]{j} =
\begin{cases}
\ecell[1]{j}                      & \text{if $j$ is odd,} \\
\ecell[1]{j} -2\ecell{(j-1,1)}    & \text{if $j>0$ is even.}
\end{cases}
\label{ecellj2undertau.eq}
\end{equation}

To find $(\expinf\tau)_\sharp\ecell[2]{j}$ using~\eqref{gensoverQ.eq}
we first need $(\expinf\tau)_\sharp\ecell[2]{1}=\ecell[1]{1}$
and $(\expinf\tau)_\sharp\ecell[2]{2}=\ecell[1]{2}-2\ecell[1]{1}\ecell[2]{1}$,
each easily found directly.  Then
\[
(\expinf\tau)_\sharp\ecell[2]{2\ell} = \frac{1}{\ell!}
  (\ecell[1]{2}-2\ecell[1]{1}\ecell[2]{1})^\ell.
\]
The cell \ecell[1]{2}\ commutes with everything, so the binomial 
theorem applies, but \ecell[1]{1}\ and \ecell[2]{1}\ square to zero, so
only two terms are nonzero. We get
\begin{align*}
(\expinf\tau)_\sharp\ecell[2]{2\ell}
    &=\frac{1}{\ell!}\left(\ecell[1]{2}\right)^\ell -
      \frac{2}{\ell!}\ell\ecell[1]{1}\ecell[2]{1}
                                       \left(\ecell[1]{2}\right)^{\ell -1} \\
    &=\ecell[1]{2\ell} - 2\ecell[1]{2\ell-1}\ecell[2]{1},
\end{align*}
the even case of~\eqref{ecellj2undertau.eq}. Multiplying by 
$(\expinf\tau)_\sharp\ecell[2]{1}=\ecell[1]{1}$ kills the second term 
and we get the odd case also.

To complete the calculation of $(\expinf\tau)_\sharp$ we find the image of
the cells $\ecell{(\ell,m)}$. If $m=2p+1$ is odd we have simply
\[
(\expinf\tau)_\sharp\ecell{(\ell,2p+1)} = \ecell[2]\ell\ecell[1]{2p+1}
    = (\expinf{\bar\tau})_\sharp\ecell{(\ell,2p+1)}, 
\]
while if $m=2p>0$ is even we get
\begin{align*}
(\expinf\tau)_\sharp\ecell{(\ell,2p)}
  &= \ecell[2]{\ell}\left(\ecell[1]{2p}-2\ecell[1]{2p-1}\ecell[2]{1}\right) \\
  &= \ecell[2]{\ell}\ecell[1]{2p}-2(-1)^\ell{\textstyle\sbinom{\ell+1}{1}}
                                       \ecell[1]{2p-1}\ecell[2]{\ell+1} \\
  &=\begin{cases}
    \ecell[2]{\ell}\ecell[1]{2p} & \text{if $\ell$ is odd}, \\
    \ecell[2]{\ell}\ecell[1]{2p}-2\ecell{(2p-1,\ell+1)} & 
                                   \text{if $\ell$ is even}.
    \end{cases}
\end{align*}
Thus
\begin{equation}
(\expinf\tau)_\sharp\ecell{(\ell,m)} =
\begin{cases}
(\expinf{\perm\tau})_\sharp\ecell{(\ell,m)}-2\ecell{(m-1,\ell+1)} 
                        & \text{$\ell,m$ both even, $m>0$}, \\
(\expinf{\perm\tau})_\sharp\ecell{(\ell,m)} & \text{otherwise}.
\end{cases}
\label{chainmapfortau.eq}
\end{equation}
Since elements of homology are linear combinations of cells each
having at least one odd index we have 
$(\Exp{k}{\tau})_*=(\Exp{k}{\perm{\tau}})_*$ as expected.

\section{The action of the braid group}
\label{braid.sec}

\subsection{Introduction}
\label{braidintro.sec}

The braid group on $n$ strands $B_n$ may be defined as the mapping class
group of an $n$--punctured disc $D_n$, or more precisely as the group of 
homeomorphisms of $D_n$ that fix $\partial D_n$ pointwise, modulo those
isotopic to the identity $\operatorname{rel} \partial D_n$. As such it acts on 
$H_k(\Exp{k}{D_n})$, and since $D_n\simeq\Gamma_n$ we may regard this
as an action on $H_k(\expg{k}{n}$). We prove the following structure 
theorem for this action.

\begin{theorem}
\label{nilpotent.th}
The image of the pure braid group $P_n$ under the action of 
$B_n$ on $H_k(\expg{k}{n})$ is nilpotent  of class at most
$\min\{(k-1)/2,\lfloor (n-1)/2\rfloor\}$ if $k$ is odd, or
$\min\{(k-2)/2,\lfloor (n-2)/2\rfloor\}$ if $k$ is even.
\end{theorem}

\noindent
For the above and other definitions of the braid group see 
Birman~\cite{birman-blmcg}.

Recall that the pure braid group $P_n$ is the kernel of the map
$B_n\rightarrow S_n$ sending each braid $\beta$ to the induced permutation
$\perm\beta$ of the punctures. Consider the subgroup
of $P_n$ consisting of braids whose first $n-1$ strands form the 
trivial braid. The $n$th strand of such a braid may be regarded as an 
element of $\pi_1(D_{n-1})$, and doing so gives an isomorphism from this
subgroup to the free group $F_{n-1}$. This shows that $P_n$ is not
nilpotent for $n\geq 3$. The group $P_2$ inside $B_2$ is isomorphic to
$2\integer$ inside \integer, and is therefore nilpotent of class $1$; 
however the bound for $n=2$ in Theorem~\ref{nilpotent.th} is zero, 
implying $P_2$ acts trivially, and in fact this follows from the second
example of section~\ref{examples.sec}. Thus we have in particular that the 
action of $B_n$ on $H_k(\expg{k}{n})$ is unfaithful for all $k$ and
$n\geq 2$. 

There is an obvious action of $S_n$ on $H_k(\expg{k}{n})$, induced
by permuting the edges of $\Gamma_n$. The theorem will be proved by
relating the action of each braid $\beta$ to that of $\perm\beta$.
Note that there is again nothing lost by considering only the action
on $H_k(\expg{k}{n})$, because of the isomorphisms induced by 
$i$ and $\cup\{v\}$. 
For brevity, in what follows we shall simply write $H_k$ for 
$H_k(\expg{k}{n})$.

\subsection{Proof of the structure theorem}
\label{nilpotent.sec}

\begin{figure}
\begin{center}
\leavevmode
\small
\psfrag{e1}{$e_1$}
\psfrag{e2}{$e_2$}
\psfrag{e3}{$e_3$}
\psfrag{e4}{$e_4$}
\psfrag{e5}{$e_5$}
\psfrag{(a)}{(a)}
\psfrag{(b)}{(b)}
\includegraphics{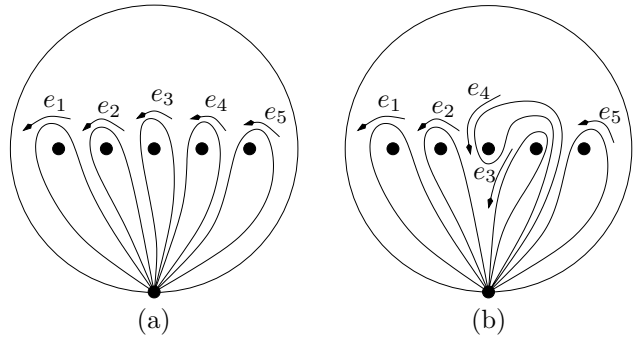}
\caption[Generators of the braid group]
{The punctured disc $D_n$ and the generators $\tau_i$ for
$n=5$ and $i=3$. We embed
$\Gamma_5$ in $D_5$ as shown in (a); the effect of $\tau_3$ on the
embedding is shown in (b).}
\label{disc.fig}
\end{center}
\end{figure}

We fix a representation of $D_n$ and a homotopy equivalence 
$\Gamma_n\rightarrow D_n$, the embedding shown in figure~\ref{disc.fig}(a).
$B_n$ is generated by the ``half Dehn twists'' $\tau_1,\ldots,\tau_{n-1}$,
where $\tau_i$ interchanges the $i$th and $(i+1)$th punctures with an
anti-clockwise twist. The effect of $\tau_i$ on the embedded graph is
shown in figure~\ref{disc.fig}(b), and we see that it induces 
the map 
\[
e_a\mapsto
\begin{cases}
e_{i+1}                 & \text{if $a=i$,} \\
e_{i+1}e_i\bar{e}_{i+1} & \text{if $a=i+1$,}\\
e_a                     &\text{if $a\not= i,i+1$,}
\end{cases}
\]
on $\Gamma_n$; regarding the $e_i$ as generators of the free group $F_n$ 
this is the standard embedding 
$B_n\hookrightarrow \operatorname{Aut}(F_n)$. We will call the induced map 
$\Gamma_n\rightarrow \Gamma_n$ $\tau_i$ also, and 
our goal is to understand the $(\Exp{k}{\tau_i})_*$ well enough
to show that the pure braid group acts by upper-triangular matrices.

To this end we define a filtration of $H_k$. Each homology class has
precisely one representative in $\tilde\C_k$, since there are no boundaries,
so we may regard $H_k$ as a subspace of $\tilde\C_k$ and work unambiguously
on the level of chains. Let
\[
\F_j = \left\{ \textstyle{\sum_J} c_J\ecell{J} \in H_k \big| 
             \text{$c_J=0$ if $|J|_2 < n-j$} \right\},
\]
so that $\F_j$ is the subspace spanned by cells having no more
than $j$ even indices. Clearly 
$\{0\}\subseteq \F_0\subseteq \F_1 \subseteq\cdots\subseteq \F_{n}=H_k$,
so the $\F_j$ form a filtration. Moreover 
$\F_j$ has $\F_j\cap \mathcal{B}(k,n)$ as a basis:
each element of $\mathcal{B}(k,n)$ has the form
\[
\widetilde{\partial\vcell{J}} = 
        \sum_{i:j_i>0\text{ even}} c_i\ecell{\drop_i(J)}
\]
for some $J$ with $|J|=k+1$, and the mod 2 norm of every term in this
sum is $|J|_2+1$. 

We now describe when $\F_i=\F_j$ for $i\not=j$, as
the indexing has been chosen to be uniform and
meaningful at the expense of certain a priori isomorphisms among the 
$\F_j$, arising from the parity and size of $k$.
Since $|J|=k$ for each $k$--cell $\ecell{J}$ we have in particular
$|J|_2 \equiv k\bmod 2$. Thus $\F_{j+1}=\F_j$ if $n-j$ and $k$ have the
same parity. Next, $|J|_2\leq |J| = k$, so $\F_j=\{0\}$ if $n-j>k$, or in
other words if $j<n-k$. Lastly, no cell for which $|J|_2 = 0$ is a summand
of an element of $H_k$ (recall that $\C_*$ is exact, and that such a cell
sits at the top of an $m$--cube complex and so is not part of any boundary),
so that either $\F_{n-2}$ or $\F_{n-1}$ is all of $H_k$, depending on the
parity of $k$. To see that these are the only circumstances in which 
$\F_i$ and $\F_j$ can coincide suppose that $i<j\leq n-1$, $n-i$ and $n-j$ 
have the same parity as $k$,  $n-j\leq k$, and let
\[
J = (\,\underbrace{0,\ldots,0}_j,\underbrace{1,\ldots,1}_{n-j-1},k+j+1-n).
\]
Every non-zero index of $J$ is odd so \ecell{J}\ is a cycle. It has
$j>i$ even indices and is therefore non-trivial in the quotient 
$\F_j/\F_{i}$.

The significance of the filtration is given by the following lemma,
which is the main step in proving the structure theorem for the action.

\begin{lemma}
Each $\F_j$ is an invariant subspace for the action of $B_n$ on $H_k$.
Moreover, the action on $\F_j/\F_{j-1}$ factors through the symmetric group
$S_n$. 
\label{filtration.lem}
\end{lemma}

\begin{proof}
The lemma is proved by relating the action of $B_n$ to that of $S_n$ 
induced by permuting the edges of $\Gamma_n$. Much of the work has
been done already in section~\ref{examples.sec}, since $\tau_i$ and
$\perm{\tau}_i$ act
as the maps $\tau$ and $\perm\tau$ considered there on $e_i\cup e_{i+1}$,
and as the identity on the remaining edges.

For each $n$--tuple $J$ let 
\[
\lift_i(J) = (j_1,\ldots,j_i+1,\ldots,j_n)
\]
so that
\[
\drop_i\circ\lift_{i+1}\circ\perm{\tau}_i(J) =
(j_1,\ldots,j_{i-1},j_{i+1}-1,j_i+1,j_{i+2},\ldots,j_n).
\]
By~\eqref{chainmapfortau.eq} we may write
\begin{equation}
(\expinf\tau_i)_\sharp\ecell{J} =
\begin{cases}
(\expinf{\perm{\tau}_i})_\sharp\ecell{J}-
       2\ecell{\drop_i\circ\lift_{i+1}\circ\perm{\tau}_i(J)} 
                        & \text{$j_i,j_{i+1}$ both even, $j_{i+1}>0$}, \\
(\expinf{\perm{\tau}_i})_\sharp\ecell{J} & \text{otherwise}.
\end{cases}
\label{braidaction.eq}
\end{equation}
If $j_i$ and $j_{i+1}$ are both even then
\[
|\drop_i\circ\lift_{i+1}\circ\perm{\tau}_i(J)|_2 = |J|_2+2
\]
and it follows that
\[
(\Exp{k}{\tau_i})_* c \in (\Exp{k}{\perm{\tau}_i})_* c + \F_{j-1}
\]
for each $c\in\F_j$ and $i\in[n-1]$. Since the $\tau_i$ generate $B_n$
we get
\[
(\Exp{k}\beta)_* c \in (\Exp{k}{\perm\beta})_* c + \F_{j-1}
\]
for all braids $\beta$ and $c\in\F_j$, and the lemma follows immediately.
\end{proof}

\begin{proof}[Proof of Theorem~\ref{nilpotent.th}]
By Lemma~\ref{filtration.lem}, with respect to a suitable ordering
of $\mathcal{B}(k,n)$ the braid group acts by block upper-triangular matrices.
The diagonal blocks are the matrices of the action on $\F_j/\F_{j-1}$, and
since this factors through $S_n$ we have that the pure braids act by
upper-triangular matrices with ones on the diagonal. It follows immediately 
that the image of $P_n$ is nilpotent.

To bound the length of the lower central series we count the number
of nontrivial blocks, as the class of the image is at most one less than this.
By the discussion following the definition of 
$\{\F_j\}$ this is the number of $0\leq j\leq n-1$ such that
$n-j\leq k$ and $n-j\equiv k\bmod 2$; letting $\ell=n-j$ this is the
number of $1\leq\ell\leq\min\{n,k\}$ such that $\ell\equiv k\bmod 2$. There
are $\lfloor (m+1)/2\rfloor$ positive odd integers and $\lfloor m/2\rfloor$
positive even integers less than or equal to a positive integer $m$, and the
given bounds follow.
\end{proof}

\subsection{The action of $B_3$ on $H_3(\expg{3}{3})$}
\label{actiononh3exp3g3.sec}

We study the action of $B_3$ on $H_3(\expg{3}{3})$, 
being the smallest non-trivial example,
and show that $P_3$ acts as a free abelian group of rank two.

For simplicity we will simply write $J$ for the cell $\ecell{J}$. From
Theorem~\ref{basis.th} we obtain the basis
\begin{align*}
u_1 &= (3,0,0) & w_1 = (0,1,2)+(0,2,1) \\
u_2 &= (0,3,0) & w_2 = (1,0,2)+(2,0,1) \\
u_3 &= (0,0,3) & w_3 = (1,2,0)+(2,1,0) \\
v   &= (1,1,1) &
\end{align*}
for $H_3(\expg{3}{3})$ (in fact this is the negative of the basis given
there). Let $U$ be the span of $\{u_1,u_2,u_3\}$, 
$V$ the span of $\{v\}$,
and $W$ the span of $\{w_1,w_2,w_3\}$. 
The subspaces $U$, $V$ and $V\oplus W$ are
easily seen to be invariant using equation~\eqref{braidaction.eq}, and
moreover the action on $U$ is simply the permutation representation 
of $S_3$. We therefore restrict our attention to $V\oplus W$. The actions 
on $V$ and $(V\oplus W)/V$ are the sign and permutation representations 
of $S_3$ respectively, and with respect to the basis $\{v,w_1,w_2,w_3\}$
we find that
\begin{align*}
(\Exp{3}{\tau_1})_*\big|_{V\oplus W} &= T_1 = 
\begin{bmatrix}
-1 & -2 & 0 & 0 \\
0  &  0 & 1 & 0 \\
0  &  1 & 0 & 0 \\
0  &  0 & 0 & 1
\end{bmatrix}, \\
(\Exp{3}{\tau_2})_*\big|_{V\oplus W} &= T_2 = 
\begin{bmatrix}
-1 &  0 &-2 & 0 \\
0  &  1 & 0 & 0 \\
0  &  0 & 0 & 1 \\
0  &  0 & 1 & 0
\end{bmatrix}.
\end{align*}
In each case the inverse is obtained by moving the $-2$ one place to the
right.

A product of $T_1$, $T_2$ and their inverses has the form
\[
P = \begin{bmatrix}
    \det\perm{P} &     p    \\
         0       & \perm{P} 
    \end{bmatrix}
\]
where $\perm{P}$ is a permutation matrix and $p=(p_1,p_2,p_3)$ is a vector of
even integers. Consider $\Sigma(P)=p_1+p_2+p_3$. Multiplying $P$ by
$T_i^{\pm 1}$ we see that $\Sigma(T_i^{\pm 1}P)=-\Sigma(P)-2$, 
and it follows that
$\Sigma(P)$ is zero if $\perm{P}$ is even and $-2$ if $\perm{P}$ is odd.
In particular $\Sigma(P)$ is zero if $P$ is the image of a pure braid.

If $P=QR$ where $\perm{Q}=\perm{R}=I$ then $p=q+r$ and it follows that
the pure braid group $P_3$ acts as a free abelian group of rank at most
two. That the rank is in fact two can be verified by calculating
$T_1^2$ and $T_2^2$ and checking that the corresponding vectors are
independent.

\appendix

\section{Table of Betti numbers}
\label{bettinumbers.apdx}

Table~\ref{bettinumbers.tab} lists the Betti numbers 
\begin{align*}
b_k(\expg{k}{n}) &=  \sum_{j=1}^k (-1)^{j-k} \binom{n+j-1}{n-1} \nonumber \\
                 &= \begin{cases}
                    \sum_{j=1}^{\ell}\binom{n+2j-2}{n-2} & 
                                              \text{if $k=2\ell$ is even}, \\
                    n+\sum_{j=1}^{\ell}\binom{n+2j-1}{n-2} &
                                        \text{if $k=2\ell+1$ is odd},
                    \end{cases}
\end{align*}
for $1\leq k\leq 20$, $1\leq n\leq 10$. To find the other non-vanishing
Betti numbers recall that
$b_{k-1}(\expg{k}{n})=b_{k-1}(\expgb{k}{n})=b_{k-1}(\expg{k-1}{n})$ 
for $k\geq 2$. 

\begin{table}[tb]
\small
\begin{tabular}{|rr|rrrrrrrrrr|} \hline
 & & \multicolumn{10}{c|}{$n$} \\
 &  & 1 & 2 & 3 & 4 & 5 & 6 & 7 & 8 & 9 & 10 \\ \hline
& 1 & 1 & 2 & 3 & 4 & 5 & 6 & 7 & 8 & 9 & 10 \\
& 2 & 0 & 1 & 3 & 6 & 10 & 15 & 21 & 28 & 36 & 45 \\
& 3 & 1 & 3 & 7 & 14 & 25 & 41 & 63 & 92 & 129 & 175 \\
& 4 & 0 & 2 & 8 & 21 & 45 & 85 & 147 & 238 & 366 & 540 \\
& 5 & 1 & 4 & 13 & 35 & 81 & 167 & 315 & 554 & 921 & 1462 \\
& 6 & 0 & 3 & 15 & 49 & 129 & 295 & 609 & 1162 & 2082 & 3543 \\
& 7 & 1 & 5 & 21 & 71 & 201 & 497 & 1107 & 2270 & 4353 & 7897 \\
& 8 & 0 & 4 & 24 & 94 & 294 & 790 & 1896 & 4165 & 8517 & 16413 \\
& 9 & 1 & 6 & 31 & 126 & 421 & 1212 & 3109 & 7275 & 15793 & 32207 \\
$k\!\!$&10& 0 & 5 & 35 & 160 & 580 & 1791 & 4899 & 12173 & 27965 & 60171 \\
&11 & 1 & 7 & 43 & 204 & 785 & 2577 & 7477 & 19651 & 47617 & 107789 \\
&12 & 0 & 6 & 48 & 251 & 1035 & 3611 & 11087 & 30737 & 78353 & 186141 \\
&13 & 1 & 8 & 57 & 309 & 1345 & 4957 & 16045 & 46783 & 125137 & 311279 \\
&14 & 0 & 7 & 63 & 371 & 1715 & 6671 & 22715 & 69497 & 194633 & 505911 \\
&15 & 1 & 9 & 73 & 445 & 2161 & 8833 & 31549 & 101047 & 295681 & 801593 \\
&16 & 0 & 8 & 80 & 524 & 2684 & 11516 & 43064 & 144110 & 439790 & 1241382 \\
&17 & 1 & 10 & 91 & 616 & 3301 & 14818 & 57883 & 201994 & 641785 & 1883168 \\
&18 & 0 & 9 & 99 & 714 & 4014 & 18831 & 76713 & 278706 & 920490 & 2803657 \\
&19 & 1 & 11 & 111 & 826 & 4841 & 23673 & 100387 & 379094 & 1299585 & 4103243\\
&20 & 0 & 10 & 120 & 945 & 5785 & 29457 & 129843 & 508936 & 1808520 & 5911762
\\ \hline
\end{tabular}
\caption{Betti numbers $b_k(\expg{k}{n})$ for $1\leq k\leq 20$ and 
         $1\leq n\leq 10$.}
\label{bettinumbers.tab}
\end{table}

\Addresses\recd

\end{document}